\documentclass[final,1p,times,11pt]{article}

\usepackage{enumerate}
\usepackage{epsfig,alltt}
\usepackage[english]{babel}
\usepackage{blindtext}
\usepackage{dsfont}
\usepackage{amssymb}
\usepackage{amsfonts}
\usepackage{graphicx}
\usepackage{amsmath}
\usepackage{amsthm}
\usepackage{natbib}
\usepackage[utf8]{inputenc}

\usepackage[left=2.5cm,right=2.5cm,top=2.5cm,bottom=2.5cm]{geometry}

\usepackage{floatflt}
\usepackage{makeidx}
\usepackage{algorithm}

\usepackage{algorithmic}

\usepackage{url}
\usepackage{float} 
\theoremstyle{remark}

\usepackage{color}

\newcommand{\begeq}[1]{\begin{equation} \label{#1}}
\newcommand{\fineq}{\end{equation}}

\newcommand{\g}{\gamma}

\title{Recursive density estimators based on Robbins-Monro's scheme and using Bernstein polynomials}

\author{Yousri SLAOUI$^{1}$\footnote{e-mail adress: Yousri.Slaoui@math.univ-poitiers.fr}~~ and ~~Asma JMAEI$^{1,2}$\footnote{e-mail adress: asma.jmaei@math.univ-poitiers.fr}\\
$^{1}$Univ. Poitiers, Lab. Math. et Appl., Futuroscope Chasseneuil, France\\
$^{2}$Facult\'e des sciences de Bizerte, Jarzouna, Tunisie}

\begin{document}

\newtheorem{theor}{Theorem}
\newtheorem{prop}{Proposition}
\newtheorem{lemma}{Lemma}
\newtheorem{lem}{Lemma}
\newtheorem{coro}{Corollary}
\newtheorem{definition}{Definition}

\newtheorem{prof}{Proof}
\newtheorem{defi}{Definition}
\newtheorem{rem}{Remark}


\date{ }
\maketitle
\begin{abstract}
\noindent \emph{
In this paper, we consider the alleviation of the boundary problem when the probability density function has bounded support. We apply Robbins-Monro's algorithm and Bernstein polynomials to construct a recursive density estimator. We study the asymptotic properties of the proposed recursive estimator. We then compared our proposed recursive estimator with many others estimators. Finally, we confirm our theoretical result through a simulation study and then using two real datasets.}

\end{abstract}

\section{Introduction}
There has been a considerable development of methods for smooth estimation of density and distribution functions, following the introduction of several kernel smoothing by \cite{Ros56} and the further advances made on kernel method by \cite{Par62}. {We advise the reader to see the paper of~\cite{{Har91}} for an introduction of several kernel smoothing techniques}. However, these methods have difficulties at and near
boundaries when curve estimation is attempted over a
region with boundaries. Moreover, it is well known in nonparametric kernel density estimation that the bias of the standard kernel density estimator 
\begin{eqnarray*}
\hat{f}(x)=\frac{1}{nh_n}\displaystyle\sum_{i=1}^nK\left(\frac{x-X_i}{h_n}\right)
\end{eqnarray*}
is of a larger order near the boundary than that in the interior, where, $K$ is a kernel (that is, a function satisfying $\int_{\mathbb{R}} K\left(x\right)dx=1$), and $\left(h_n\right)$ is a bandwidth  (that is, a sequence of positive real numbers that goes to zero). We suppose for simplicity that there is a single known boundary to the support of the density function $f$ which we might as well take to be at the origin, {then} we are dealing
with positive data. For convenience, we consider a symmetric kernel (for
instance, normal kernels).
Away from the boundary, which means that at any $x>h_n$, the usual asymptotic
mean and variance expressions are applied. Let us now suppose that $f$ has two
continuous derivatives everywhere, and that as $n\to\infty$, 
$h = h_n\to 0$ and $nh\to 0$. Then, 
\begin{eqnarray*}
\mathbb{E}\left[\hat{f}(x)\right]\simeq f(x)+\frac{1}{2}h^2f''(x)\int x^2K(x)dx,
\end{eqnarray*}
and 
\begin{eqnarray*}
Var\left[\hat{f}(x)\right]\simeq(nh)^{-1}f(x) \int K^2(x)dx.
\end{eqnarray*}
Near the boundary, the expression of the mean and the variance are different.
Let $x=ph$, we have 
\begin{eqnarray*}
\mathbb{E}\left[\hat{f}(x)\right]\simeq \displaystyle f(x)\int_{-\infty}^pK(x)dx-f'\left(x\right)\displaystyle\int_{-\infty}^pxK(x)dx+\frac{1}{2}h^2\displaystyle f''(x)\int_{-\infty}^px^2K(x)dx,
\end{eqnarray*}
and 
\begin{eqnarray*}
Var\left[\hat{f}(x)\right]\simeq(nh)^{-1}\displaystyle f(x)\int_{-\infty}^pK^2(x)dx.
\end{eqnarray*}
These bias phenomena are called boundary bias. Many authors have suggested methods for reducing this phenomena such as data reflection (\cite{Sch85}), boundary kernels (\cite{Mul91, Mul93} and \cite{Mul94}), the local linear estimator (\cite{Lej92} and \cite{Jon93}), the use of beta and gamma kernels (\cite{Che99, Che00}).\\
For a smooth estimate of a density function with finite known support, Vitale's method (\cite{Vit75}) based on the Bernstein polynomials, illustrated below, also has been investigated in the literature (\cite{Gho01}, \cite{Bab02}, \cite{Kak04}, \cite{Rao05}) and more recently by~\cite{Leb10} and \cite{Kak14}). The idea comes from the Weierstrass's approximation theorem that for any continuous function $u$ on the interval $[0,1]$, we have 
\begin{eqnarray*}
\displaystyle\sum_{k=0}^m u\left(\frac{k}{m}\right)b_k(m,x)\to u(x),\quad \mbox{uniformly in}\, x\in[0,1], 
\end{eqnarray*}
where $b_{k}(m,x)={{m}\choose{k}}x^{k}(1-x)^{m-k}$ is the Bernstein polynomial of order $m$ .\\
In the context of distribution function $F$ with support $[0,1]$, \cite{Vit75} proposed an estimator 
\begin{eqnarray*}
\widetilde{F}_n(x)=\displaystyle\sum_{k=0}^{m}F_n\left(\frac{k}{m}\right)b_k(m,x),
\end{eqnarray*}
where $F_n$ is the empirical distribution based on a random sample $X_{1},X_{2}, \ldots, X_{n}$. Hence, an estimator for the density $f$ is given by
\begin{eqnarray}
\widetilde{f}_n(x)&=&\frac{d}{dx}\widetilde{F}_n(x)=m\displaystyle\sum_{k=0}^{m-1}\left[F_n\left(\frac{k+1}{m}\right)-F_n\left(\frac{k}{m}\right)\right]b_k(m-1,x).
\label{eq:26}
\end{eqnarray}
In this paper, we propose a recursive method to estimate the unknown density function $f$. The advantage of recursive estimators is that their update from a sample of size $n$ to one of size
$n + 1$, requires considerably less computations. This property is particularly important, since the number of points at which the function is estimated is usually very large.\\
$ $\\
Let $X_{1},X_{2}, \ldots, X_{n}$ be a sequence of i.i.d random variables having a common unknown distribution $F$ with associated density $f$ supported on $[0,1]$. {In order to construct a recursive method to estimate the unknown density $f$, we first follow \cite{Leb10} and we introduce $T_{n,m}$ as follows}: 
\begin{eqnarray*}
T_{n,m}(x)&=&m\sum^{m-1}_{k=0}\left(\mathbb{I}_{\left\{X_n\leq\frac{k+1}{m}\right\}}-\mathbb{I}_{\left\{X_{n}\leq\frac{k}{m}\right\}}\right)b_{k}(m-1,x)\\
&=&m\sum^{m-1}_{k=0}\mathbb{I}\left\{\frac{k}{m}<X_{n}\leq\frac{k+1}{m}\right\}b_{k}(m-1,x)\\
&=&mb_{k_n}(m-1,x).
\end{eqnarray*}
Here, we let $k_n=[mX_n]$, where $[x]$ denotes the largest integer smaller than $x$.
Then, we use Robbins-Monro's scheme (\cite{Rob51}), and we set $f_{0}(x)\in\mathbb{R}$ and {for all $n\in\mathbb{N}^*$, we set} 
\begin{eqnarray}\label{eq:rec:density}
f_{n}(x)=(1-\g_n)f_{n-1}(x)+\gamma_{n}Z_{n}\left(x\right),
\end{eqnarray}
{where $(\g_n)$ is a sequence of real numbers, called a stepsize and $Z_{n}(x)=2T_{n,m}(x)-T_{n,m/2}(x)$. Then for simplicity, we suppose that $f_0(x)=0$ and $\Pi_n=\prod_{j=1}^n(1-\gamma_j)$}.

Then, it follows from~\eqref{eq:rec:density}, that one can estimate $f$ recursively at the point $x$ by 
\begin{eqnarray}
f_n(x)=\Pi_n\sum_{k=1}^n\Pi_k^{-1}\gamma_k Z_k(x).
\label{eq:1}
\end{eqnarray}
Our first aim in this paper is to compute the bias, the variance, the mean squared error ($MSE$) and the mean integrated squared error ($MISE$) of our proposed recursive estimators. It turns out that they heavily depend on the choice of the stepsize $(\g_n)$. Moreover, we give the optimal order $(m_n)$ which minimizes the $MSE$ and the $MISE$ of the proposed recursive estimators. Further, we show that using the stepsize $(\g_n)=(n^{-1})$ and the optimal order $\left(m_n\right)$, the proposed estimator $f_n$ can dominate Vitale's estimator $\widetilde{f}_n$ in terms of $MISE$. Finally, we confirm our theoretical results from a simulation study.\\

The remainder of this paper is organized as follows. In Section~\ref{section:notation:assumption}, we list our notations and assumptions. In Section~\ref{section:main:result}, we state the main theoretical results regarding bias, variance, $MSE$ and $MISE$. Section~\ref{section:app} is devoted to some numerical studies : first by simulation (Subsection~\ref{subsection:simu}) and second using some real datasets (Subsection~\ref{subsection:realdata}). We conclude the paper in Section~\ref{section:conclusion}. Appendix~\ref{section:proof} gives the proof of our theoretical results. 
\section{Assumptions and Notations} \label{section:notation:assumption}
\begin{defi}\label{def:1}
Let $\gamma\in\mathbb{R}$ and $(v_n)_{n\ge 1}$ be a nonrandom positive sequence. We say that $(v_n)\in \mathcal{GS}(\gamma)$ if \\
\begin{equation}
\lim_{n \rightarrow +\infty} n\left[1-\frac{v_{n-1}}{{v_n}}\right]=\gamma. \nonumber
\end{equation}
\end{defi}
This condition was introduced by \cite{Gal73} to define regularly varying sequences (see also \cite{Boj73}). Typical sequences in $\mathcal{GS}(\gamma)$ are, for $b\in \mathbb{R}$, $n^\gamma(\log n)^b$, $n^\gamma(\log \log n)^b$, and so on.\\
To study the estimator (\ref{eq:1}), we make the following assumptions :
\begin{description}
\item $\left(A1\right)$ $f$ admits a continuous fourth-order derivative $f^{(4)}$ on $[0, 1]$.
\item $\left(A2\right)$ $(\gamma_{n})\in \mathcal{GS}\left(-\alpha\right)$, $\alpha\in(\frac{1}{2},1]$.
\item $\left(A3\right)$ $(m_{n})\in \mathcal{GS}(a)$, $a\in(0,1)$.
\item$\left(A4\right)$ $\lim_{n\rightarrow\infty}(n\gamma_{n})\in(\min\left(2a,(2\alpha-a)/4\right),\infty)$.
\end{description}
\begin{itemize}
\item Assumption $(A1)$ is standard in the framework of nonparametric estimation of probability density using Bernstein polynomials (see \cite{Leb10}).
\item Assumption $(A2)$ on the stepsize is usual in the framework of the recursive estimation for density estimation (see \cite{Mok09} and \cite{Sla13,Sla14a,Sla14b}). {This assumption ensures that $\sum_{n\geq 1} \gamma_n=\infty$  and $\sum_{n\geq 1} \gamma_n^2<\infty$, which are two classical assumptions for obtaining the convergence of Robbins-Monro's algorithm  (see~\cite{Duf97})}. 
\item Assumption $(A3)$ on $(m_n)$ was introduced similarly to the assumption on the bandwidth used for the recursive kernel distribution estimator (see \cite{Sla14a,Sla14b}) {to ensure the application of the technical lemma given in the appendix A}. 
\item Assumption $(A4)$ on the limit of $\left(n\gamma_n\right)$ as $n$ goes to infinity is usual in the framework of stochastic approximation algorithms. {This condition ensures the application of the technical lemma given in the appendix A to obtain the asymptotic bias and variance respectively}.
\end{itemize}
\begin{rem}
The intuition behind the use of such order $\left(m_n\right)$ belonging to $\mathcal{GS}\left(a\right)$ is that the ratio $m_{n-1}/m_n$ is equal to $1-a/n+o\left(1/n\right)$, then using such order and using the assumption on the stepsize, which ensures that  $\gamma_{n-1}/\gamma_n$ is equal to $1+\alpha/n+o\left(1/n\right)$, the application of the technical lemma given in the appendix A ensures that the bias and the variance will depend only on $m_n$ and $\gamma_n$ and not on $m_1,\ldots,m_n$ and $\gamma_1,\ldots,\gamma_n$, then the $MISE$ will depend also only on $m_n$ and $\gamma_n$, which will be helpful to deduce an optimal order and an optimal stepsize.
\end{rem}
Throughout this paper we will use the following notations:
\begin{eqnarray*}
\Delta_1(x)&=&\frac{1}{2}\left[(1-2x)f'(x)+x(1-x)f''(x)\right],\quad \psi(x)=(4\pi x(1-x))^{-1/2}, \quad \xi=\displaystyle\lim_{n\to\infty}(n\g_n)^{-1},\\
\Delta_2(x)&=&\frac{1}{6}(1-6x(1-x))f''(x)+\frac{5}{12}x(1-x)(1-2x)f^{(3)}(x)+\frac{1}{8}x^2(1-x)^2f^{(4)}(x),\\
C_1&=&\int_0^1f(x)\psi(x)(x)dx,\quad C_2=\int_0^1\Delta_2^2(x)dx, \quad C_3=\frac{1}{\sqrt{2}}+4\left(1-\sqrt{\frac{2}{3}}\right),\\ C_4&=&\int_0^1\Delta_1^2(x)dx, \quad C_5=\displaystyle\int_0^1\left\{-\frac{\Delta^2_1(x)}{2f(x)}+\Delta_2(x)\right\}^2dx,\\ C_6&=&\displaystyle\int_0^1\left\{-\frac{\Delta^2_1(x)}{2f(x)}+\Delta_2(x)+f(x)\displaystyle\int_0^1\frac{\Delta^2_1(y)}{2f(y)}dy\right\}^2dx.
\end{eqnarray*}
\section{Main results} \label{section:main:result}
{In this section we start} by developing the bias, the variance and the $MSE$ of the proposed recursive estimators using Bernstein estimator first in the boundary region and then at the edges.
\subsection{Within the interval $[0,1]$}\label{subsection:with}
The following proposition gives the bias, variance, and $MSE$ of the proposed recursive estimator $f_n(x)$ for $x\in(0, 1)$.
\begin{prop}\label{prop:1}
$ $\\
Let Assumptions $\left(A1\right)-\left(A4\right)$ hold. For $x\in(0,1)$, we have
\begin{enumerate}
\item If $a\in\left(0,\frac{2}{9}\alpha\right]$, then
\begin{equation}
Bias\left[f_n(x)\right]=-m^{-2}_{n} \frac{2}{1-2a\xi}\Delta_2(x)+o\left(m^{-2}_{n}\right).
\label{eq:2}
\end{equation}
If $a\in\left(\frac{2}{9}\alpha,1\right)$, then
\begin{equation}
Bias\left[f_n(x)\right]=o\left(\sqrt{\gamma_{n}m^{1/2}_{n}}\right).
\label{eq:3}
\end{equation}
\item If $a\in\left[\frac{2}{9}\alpha,1\right)$, then
\begin{equation}
Var[f_{n}(x)]=C_3\gamma_{n}m^{1/2}_{n}\frac{2}{4-(2\alpha-a)\xi}f(x)\psi(x)+o\left(\gamma_{n}m^{1/2}_{n}\right).
\label{eq:4}
\end{equation}
If $a\in\left(0,\frac{2}{9}\alpha\right]$, then
\begin{equation}
Var[f_{n}(x)]=o\left(m_n^{-4}\right).
\label{eq:5}
\end{equation}
\item If $\lim_{n\rightarrow\infty}(n\gamma_{n})>\max\left(2a,(2\alpha-a)/4\right)$, that implies $\lim_{n\rightarrow\infty}(n\gamma_{n})>2a$ and $\lim_{n\rightarrow\infty}(n\gamma_{n})>(2\alpha-a)/4$, {in this case, we have}  $1-2a\xi>0$ and $4-(2\alpha-a)\xi>0$ then (\ref{eq:2}) and (\ref{eq:4}) hold simultaneously.
\item If $a\in\left(0,\frac{2}{9}\alpha\right)$, then
\begin{eqnarray*}
MSE\left[f_{n}(x)\right]=\Delta_2^2(x)m^{-4}_{n}\frac{4}{(1-2a\xi)^2}+o\left(m^{-4}_{n}\right).
\end{eqnarray*}
If $a=\frac{2}{9}\alpha$, then
\begin{eqnarray}\label{eq:MSE:2.9}
MSE\left[f_{n}(x)\right]&=&\Delta_2^2(x)m^{-4}_{n}\frac{4}{(1-2a\xi)^2}+C_3f(x)\psi(x)\gamma_{n}m^{1/2}_{n}\frac{2}{4-(2\alpha-a)\xi}\nonumber\\
&&+o\left(m^{-4}_{n}+\gamma_{n}m^{1/2}_{n}\right).
\end{eqnarray}
If $a\in\left(\frac{2}{9}\alpha,1\right]$, then
\begin{eqnarray*}
MSE\left[f_{n}(x)\right]&=&C_3f(x)\psi(x)\gamma_{n}m^{1/2}_{n}\frac{2}{4-(2\alpha-a)\xi}+o\left(\gamma_{n}m^{1/2}_{n}\right).
\end{eqnarray*}
\end{enumerate}
\end{prop}
{In order to choose the optimal order $\left(m_n\right)$, we choose $\left(m_n\right)$ such that}
\begin{eqnarray*}
{m_n=\underset{m_n}{\arg\min}MSE\left[f_{n}\left(x\right)\right],\quad \mbox{for}\,\, x\in(0,1).}
\end{eqnarray*}
{Then, it follows from~\eqref{eq:MSE:2.9}, that $\left(m_n\right)$ must be equal to} 
\begin{eqnarray*}
\left(2^{2/9}\left(1-\frac{4}{9}\xi\right)^{-2/9}\left[\frac{32\Delta_2^2(x)}{C_3f(x)\psi(x)}\right]^{2/9}\gamma_n^{-2/9}\right),
\end{eqnarray*}
{and then the corresponding $MSE$ is equal to}
\begin{eqnarray*}
MSE\left[f_n(x)\right]=\frac{9(32C_3^8)^{1/9}}{8}\frac{(\Delta_2(x))^{2/9}(f(x)\psi(x))^{8/9}}{2^{2/9}\left(1-\frac{4}{9}\xi\right)^{10/9}}\gamma_n^\frac{8}{9}+o\left(\gamma_n^\frac{8}{9}\right).
\end{eqnarray*}
{Moreover, since the optimal stepsize should be obtained by minimizing the $MSE$ of our proposed estimators $f_n\left(x\right)$, then $\left(\gamma_n\right)$ must be chosen in $\mathcal{GS}\left(-1\right)$}. {By considering the case when $(\g_n)=\left(\g_0n^{-1}\right)$, we obtain} 
\begin{eqnarray*}
\left(m_n\right)=\left(2^{2/9}\left(\g_0-4/9\right)^{-2/9}\left[\frac{32\Delta_2^2(x)}{C_3f(x)\psi(x)}\right]^{2/9}n^{2/9}\right),
\end{eqnarray*}
{and the corresponding $MSE$ is equal to}
\begin{eqnarray*}
MSE\left[f_n(x)\right]=\frac{9(32C_3^8)^{1/9}(\Delta_2(x))^{2/9}(f(x)\psi(x))^{8/9}}{8}\frac{\g_0^2}{2^{8/9}\left(\g_0-4/9\right)^{10/9}}n^{-8/9}+o\left(n^{-8/9}\right).
\end{eqnarray*}
\subsection{The edges of the interval $[0,1]$}\label{subsection:edgs}
For the cases $x=0,1$, we need {the following} additional assumption :
\begin{description}
\item $\left(A'4\right)$ $\quad\lim_{n\rightarrow\infty}(n\gamma_{n})\in(\min\left(2a,(\alpha-a)/2\right),\infty)$.
\end{description}
The following proposition gives bias, variance and $MSE$ of $f_{n}(x)$, for $x=0,1$.
\begin{prop}\label{prop:2}
$ $\\
Let Assumptions $\left(A1\right)-\left(A3\right)$ and $\left(A^{\prime}4\right)$ hold, for $x=0,1$, we have 
\begin{enumerate}
\item If $a\in\left(0,\frac{\alpha}{5}\right]$, then
\begin{equation}
Bias\left[f_n(x)\right]=-m^{-2}_{n} \frac{2}{1-2a\xi}\Delta_2(x)+o\left(m^{-2}_{n}\right).
\label{eq:2'}
\end{equation}
If $a\in\left(\frac{\alpha}{5},1\right)$, then
\begin{equation}
Bias\left[f_n(x)\right]=o\left(\sqrt{\gamma_{n}m_{n}}\right).
\label{eq:3'}
\end{equation}
\item If $a\in\left[\frac{\alpha}{5},1\right)$, then
\begin{equation}
Var[f_{n}(x)]=\frac{5}{2}\gamma_{n}m_{n}\frac{1}{2-(\alpha-a)\xi}f(x)+o\left(\gamma_{n}m_{n}\right).
\label{eq:4'}
\end{equation}
If $a\in\left(0,\frac{\alpha}{5}\right)$, then
\begin{equation}
Var[f_{n}(x)]=o\left(m_n^{-4}\right).
\label{eq:5'}
\end{equation}
\item If $\lim_{n\rightarrow\infty}(n\gamma_{n})>\max\left(2a,(\alpha-a)/2\right)$, that implies $\lim_{n\rightarrow\infty}(n\gamma_{n})>2a$ and  $\lim_{n\rightarrow\infty}(n\gamma_{n})>(\alpha-a)/2$ {then} $1-2a\xi>0$ and $2-(\alpha-a)\xi>0$, then, (\ref{eq:2'}) and (\ref{eq:4'}) hold simultaneously.
\item If $a\in\left(0,\frac{\alpha}{5}\right)$, then
\begin{eqnarray*}
MSE\left[f_{n}(x)\right]=\Delta_2^2(x)m^{-4}_{n}\frac{4}{(1-2a\xi)^2}+o\left(m^{-4}_{n}\right).
\end{eqnarray*}
If $a=\frac{\alpha}{5}$, then
\begin{eqnarray}\label{eq:MSE:1.5}
MSE\left[f_{n}(x)\right]&=&\Delta_2^2(x)m^{-4}_{n}\frac{4}{(1-2a\xi)^2}+\frac{5}{2}f(x)\gamma_{n}m_{n}\frac{1}{2-(\alpha-a)\xi}\nonumber\\
&&+o\left(m^{-4}_{n}+\gamma_{n}m_{n}\right).
\end{eqnarray}
If $a\in\left(\frac{\alpha}{5},1\right)$, then
\begin{eqnarray*}
MSE\left[f_{n}(x)\right]&=&\frac{5}{2}f(x)\gamma_{n}m_{n}\frac{1}{2-(\alpha-a)\xi}+o\left(\gamma_{n}m_{n}\right).
\end{eqnarray*}
\end{enumerate}
\end{prop}
{In order to choose the optimal order $\left(m_n\right)$, we choose $\left(m_n\right)$ such that}
\begin{eqnarray*}
{m_n=\underset{m_n}{\arg\min}MSE\left[f_{n}\left(x\right)\right],\quad \mbox{for}\,\, x=0,1}
\end{eqnarray*}
{Then, it follows from~\eqref{eq:MSE:1.5}, that $\left(m_n\right)$ must be equal to} 
\begin{eqnarray*}
\left(2^{1/5}\left(1-\frac{2}{5}\xi\right)^{-1/5}\left[\frac{32\Delta_2^2(x)}{5f(x)}\right]^{1/5}\g_n^{-1/5}\right),
\end{eqnarray*}
{and then the corresponding $MSE$ is equal to}
\begin{eqnarray*}
MSE\left[f_n(x)\right]=\frac{5^{8/5}32^{1/5}}{8}\frac{\left(\Delta_2(x)\right)^{2/5}\left(f(x)\right)^{4/5}}{2^{4/5}\left(1-\frac{2}{5}\xi\right)^{6/5}}\g_n^{4/5}+o\left(\g_n^{4/5}\right).
\end{eqnarray*}
{Moreover, since the optimal stepsize should be obtained by minimizing the $MSE$ of our proposed estimators $f_n\left(x\right)$, then $\left(\gamma_n\right)$ must be chosen in $\mathcal{GS}\left(-1\right)$}. {By considering the case when $(\g_n)=\left(\g_0n^{-1}\right)$, we obtain} 
\begin{eqnarray*}
\left(2^{1/5}\left(\g_0-2/5\right)^{-1/5}\left[\frac{32\Delta_2^2(x)}{5f(x)}\right]^{1/5}n^{1/5}\right),
\end{eqnarray*}
{and then the corresponding $MSE$ is equal to}
\begin{eqnarray*}
MSE\left[f_n(x)\right]=\frac{5^{8/5}32^{1/5}\left(\Delta_2(x)\right)^{2/5}\left(f(x)\right)^{4/5}}{8}\frac{\g_0^2}{2^{4/5}\left(\g_0-2/5\right)^{6/5}}n^{-4/5}+o\left(n^{-4/5}\right).
\end{eqnarray*}
{In the next subsection we give the $MISE$ of the proposed recursive estimator $f_n$ introduced in~\eqref{eq:rec:density} for $x\in\left(0,1\right)$.}
\subsection{$MISE$ of the recursive estimator $f_n$}\label{section:MISE}
The following proposition gives the $MISE$ of the proposed recursive estimator $f_n$.
\begin{prop}\label{prop:3}
$ $\\
Let Assumptions $\left(A1\right)-\left(A4\right)$ hold. We have 
\begin{enumerate}
\item If $a\in\left(0,\frac{2}{9}\alpha\right)$, then
\begin{eqnarray*}
MISE\left[f_{n}\right]=C_2m^{-4}_{n}\frac{4}{(1-2a\xi)^2}+o\left(m^{-4}_{n}\right).
\end{eqnarray*}
\item If $a=\frac{2}{9}\alpha$, then
\begin{eqnarray*}
MISE\left[f_{n}\right]&=&C_2m^{-4}_{n}\frac{4}{(1-2a\xi)^2}+C_1C_3\gamma_{n}m^{1/2}_{n}\frac{2}{4-(2\alpha-a)\xi}\\
&&+o\left(m^{-4}_{n}+\gamma_{n}m^{1/2}_{n}\right).
\end{eqnarray*}
\item If $a\in\left(\frac{2}{9}\alpha,1\right)$, then
\begin{eqnarray*}
MISE\left[f_{n}\right]&=&C_1C_3\gamma_{n}m^{1/2}_{n}\frac{2}{4-(2\alpha-a)\xi}+o\left(\gamma_{n}m^{1/2}_{n}\right).
\end{eqnarray*}
\end{enumerate}
\end{prop}
The following result is a consequence of the previous proposition which gives the optimal order $\left(m_n\right)$ of the estimator (\ref{eq:1}) and the corresponding $MISE$.
\begin{coro}\label{coro:1}
$ $\\
Let Assumptions $\left(A1\right)-\left(A4\right)$ hold. To minimize the $MISE$ of $f_n$, the stepsize $(\gamma_n)$ must be chosen in $\mathcal{GS}(-1)$ and $(m_n)$ must be in $\mathcal{GS}(2/9)$ such that 
\begin{eqnarray*}
\left(2^{2/9}\left(1-\frac{4}{9}\xi\right)^{-2/9}\left[\frac{32C_2}{C_1C_3}\right]^{2/9}\gamma_n^{-2/9}\right),
\end{eqnarray*}
{and then the corresponding $MISE$ is equal to} 
\begin{eqnarray*}
MISE\left[f_{n}\right]=\frac{9(32C_1^8C_3^8C_2)^{1/9}}{8}\frac{1}{2^{8/9}\left(1-\frac{4}{9}\xi\right)^{10/9}}\gamma_n^\frac{8}{9}+o\left(\gamma_n^\frac{8}{9}\right).
\end{eqnarray*}
Moreover, in the case when $(\g_n)=\left(\g_0n^{-1}\right)$, {the optimal order $\left(m_n\right)$ must be equal to}
\begin{eqnarray}
\left(2^{2/9}\left(\g_0-4/9\right)^{-2/9}\left[\frac{32C_2}{C_1C_3}\right]^{2/9}n^{2/9}\right),
\label{eq:30}
\end{eqnarray}
{and the corresponding $MISE$ is given by}
\begin{eqnarray}
MISE\left[f_{n}\right]=\frac{9\left(8C_1^8C_3^8C_2\right)^{1/9}}{8}\frac{\g_0^2}{2^{6/9}\left(\g_0-4/9\right)^{10/9}}n^{-8/9}+o\left(n^{-8/9}\right).
\end{eqnarray}
\label{eq:31}
\end{coro}
\begin{rem}
The minimum of $\frac{\g_0^2}{\left(\g_0-4/9\right)^{10/9}}$ is reached for $\g_0=1$. Moreover, to minimize the variance of $f_n$, we should choose $\g_0=1-\frac{a}{2}$, with $a=2/9$ for $x\in(0,1)$ and $\g_0=1-a$ with $a=1/5$ for $x=0,1$. Therefore, in {the application section}, we will consider the following stepsizes $(\g_n)=(n^{-1})$, $(\g_n)=\left(\frac{8}{9}n^{-1}\right)$ and $(\g_n)=\left(\frac{4}{5}n^{-1}\right)$.
\end{rem}
Let us now state the following theorem, which gives the weak convergence rate of the estimator $f_n$ defined in (\ref{eq:1}).
\begin{theor}[Weak pointwise convergence rate]\label{theo:1}
$ $\\
Let Assumption $\left(A1\right)-\left(A4\right)$ hold. For $x\in(0,1)$, we have \\
\begin{enumerate}
\item If $\gamma^{-1/2}_{n}m^{-9/4}_{n}\rightarrow c$ for some constant $c\geq 0$, then \\
\begin{eqnarray*}
\gamma^{-1/2}_{n}m^{-1/4}_{n}(f_{n}(x)-f(x))\stackrel{\mathcal{D}}{\rightarrow}\mathcal{N}\left(-\frac{2c}{1-2a\xi}\Delta_2(x),\frac{2}{4-(2\alpha-a)\xi}C_3f(x)\psi(x)\right).
\end{eqnarray*}
\item If $\gamma^{-1/2}_{n} m^{-9/4}_{n}\rightarrow \infty $, then \\
\begin{eqnarray*}
m_{n}^{-2}\left(f_{n}(x)-f(x)\right)\stackrel{\mathbb{P}}{\rightarrow}-\frac{2}{1-2a\xi}\Delta_2(x),
\end{eqnarray*}
\end{enumerate}
where $\stackrel{\mathcal{D}}{\rightarrow}$ denotes the convergence in distribution, $\mathcal{N}$ the Gaussian-distribution and $\stackrel{\mathbb{P}}{\rightarrow}$ the convergence in probability.
\end{theor}
\begin{rem}
When the bandwidth $\left(h_n\right)$ is chosen such that $\lim_{n\to \infty}\gamma^{-1/2}_{n}m^{-9/4}_{n}=0$ and the stepsize such that $\lim_{n\to \infty}n\gamma_n=\gamma_0$, the proposed recursive estimators fulfills the following central limit theorem
\begin{eqnarray*}
n^{1/2}m^{-1/4}_{n}(f_{n}(x)-f(x))\stackrel{\mathcal{D}}{\rightarrow}\mathcal{N}\left(0,\frac{\gamma_0^2}{2\gamma_0-8/9}C_3f(x)\psi(x)\right).
\end{eqnarray*} 
\end{rem}

\subsection{Results on some classical Bernstein density estimator}
{In the next paragraph, we recall some results on Vitale's Bernstein density estimator $\widetilde{f}_n$ given in~\eqref{eq:26}}. 
\subsubsection{Vitale's Bernstein density estimator $\widetilde{f}_n$}
Under some classical assumption on the density, such as $f$ is continuous and admits two continuous and bounded derivatives, for $x\in[0,1]$, we have 
\begin{equation*}
Bias\left[\widetilde{f}_n(x)\right]=\frac{\Delta_1(x)}{m}+o\left(m^{-1}\right), \text{uniformly in}\, x\in[0,1].
\end{equation*}
\begin{equation*}
Var[\widetilde{f}_{n}(x)]=\begin{cases}
\frac{m^{1/2}}{n}f(x)\psi(x)+o_x\left(\frac{m^{1/2}}{n}\right)& \quad\text{for }x\in(0,1),\\
\frac{m}{n}f(x)+o_x\left(\frac {m}{n}\right)&\quad\text{for }x=0, 1.
\end{cases}
\end{equation*}
Moreover, we have
\begin{eqnarray*}
MISE\left[\widetilde{f}_{n}\right]&=&\frac{m^{1/2}C_1}{n}+\frac{C_4}{m^2}+o\left(m^{1/2}{n}^{-1}\right)+o\left(m^{-2}\right).
\end{eqnarray*}
To minimize the $MISE$ of $\widetilde{f}_n$, the order $(m_n)$ must equal to
\begin{eqnarray}
\left(\left[\frac{4C_4}{C_1}\right]^{2/5}n^{2/5}\right),
\label{eq:27}
\end{eqnarray}
and the corresponding $MISE$
\begin{eqnarray*}
MISE\left[\widetilde{f}_{n}\right]=\frac{5\left(C_1^4C_4\right)^{1/5}}{4}n^{-4/5}+o\left(n^{-4/5}\right).
\label{eq:28}
\end{eqnarray*}
\begin{rem}
If we let $h=m^{-1}$ be the bandwidth of the Vitale's estimator $\widetilde{f}_n$, it is clear that the bias of $\widetilde{f}_n$ is $O\left(h^{-1}\right)$ which is larger than the classical kernel density estimators which is $O(h^2)$ except near the boundaries. To reduce this bias, \cite{Leb10} suggest a new Bernstein estimator using the method of bias correction. This methodology was used first by \cite{Pol95} in the context of spectral density estimation and is linked with the work of \cite{Sch71} and \cite{Sch77} on bias reduction in estimation.
\end{rem} 
\subsubsection{Bias correction for Bernstein density $\widetilde{f}_{n,m,m/2}$}
{The bias correction for Bernstein density proposed by Leblanc~\cite{Leb10} is defined by} 
\begin{eqnarray}
\widetilde{f}_{n,m,m/2}(x)=2\widetilde{f}_{n,m}(x)-\widetilde{f}_{n,m/2}(x), \quad x\in[0,1]
\label{eq:32}
\end{eqnarray}
where $\widetilde{f}_{n,m}$ and $\widetilde{f}_{n,m/2}$ are the Bernstein density estimators introduced by Vitale with order $m$ and $m/2$ respectively which defined in (\ref{eq:26}). Let us now recall the characteristics of the estimator $\widetilde{f}_{n,m,m/2}$. Under the Assumption $(A1)$, we have 
\begin{equation*}
Bias[\widetilde{f}_{n,m,m/2}(x)]=-2\frac{\Delta_2(x)}{m^2}+o\left(m^{-2}\right),\text{ uniformly in}\, x\in[0,1].
\end{equation*}
\begin{equation*}
Var[\widetilde{f}_{n,m,m/2}(x)]=\begin{cases}
C_3\frac{m^{1/2}}{n}f(x)\psi(x)+o_x\left(\frac{m^{1/2}}{n}\right)& \quad\text{for }x\in(0,1),\\
\frac{5}{2}\frac{m}{n}f(x)+o_x\left(\frac {m}{n}\right)&\quad\text{for }x=0, 1.
\end{cases}
\end{equation*}
and then
\begin{eqnarray*}
MISE\left[\widetilde{f}_{n,m,m/2}\right]=\frac{C_1C_3m^{1/2}}{n}+\frac{4C_2}{m^4}+o\left(m^{1/2}{n}^{-1}\right)+o\left(m^{-4}\right).
\end{eqnarray*}
To minimize the $MISE$ of $\widetilde{f}_{n,m,m/2}$, the order $(m_n)$ must equal to
\begin{eqnarray}
\left(\left[\frac{32C_2}{C_1C_3}\right]^{2/9}n^{2/9}\right),
\label{eq:33}
\end{eqnarray}
{and then the corresponding $MISE$}
\begin{eqnarray*}
MISE\left[\widetilde{f}_{n,m,m/2}\right]=\frac{9\left(32C_1^8C_2C_3^8\right)^{1/9}}{8}n^{-8/9}+o\left(n^{-8/9}\right).
\label{eq:34}
\end{eqnarray*}
\cite{Kak14} have generalized the estimator proposed by Leblanc $\widetilde{f}_{n,m,m/2}$ and defined in (\ref{eq:32}) 
\begin{eqnarray}
\widetilde{f}_{n,m,m/b}(x)=\frac{b}{b-1}\widetilde{f}_{n,m}(x)-\frac{1}{b-1}\widetilde{f}_{n,m/b}(x), \quad x\in[0,1]
\label{eq:35}
\end{eqnarray}
where $b=2,3, \ldots$ and $\widetilde{f}_{n,m}$ and $\widetilde{f}_{n,m/b}$ are the Vitale's density estimators defined in (\ref{eq:26}) with order $m$ and $m/b$ respectively. Under the Assumption $\left(A1\right)$, we have  
\begin{equation*}
Bias[\widetilde{f}_{n,m,m/b}(x)]-f(x)=-\frac{b}{m^2}\Delta_2(x)+o\left(m^{-2}\right), \text{ uniformly in}\, x\in[0,1].
\end{equation*}
\begin{equation*}
Var[\widetilde{f}_{n,m,m/b}(x)]=\begin{cases}
\lambda_1(b)\frac{m^{1/2}}{n}f(x)\psi(x)+o_x\left(\frac{m^{1/2}}{n}\right)& \quad\text{for }x\in(0,1),\\
\lambda_2(b)\frac{m}{n}f(x)+o_x\left(\frac {m}{n}\right)&\quad\text{for }x=0, 1,
\end{cases}
\end{equation*} 
\begin{eqnarray*}
\lambda_1(b)&=&\frac{1}{(1-b)^2}\left\{b^2+b^{-1/2}-2b\left(\frac{2}{b+1}\right)^{1/2}\right\},\\
\lambda_2(b)&=&\frac{1}{(1-b)^2}\left\{b^2+b^{-1}-2\right\}.
\end{eqnarray*}
Moreover, we have
\begin{eqnarray*}
MISE\left[\widetilde{f}_{n,m,m/b}\right]=\lambda_1(b)C_1\frac{m^{1/2}}{n}+\frac{b^2C_2}{m^4}+o\left(m^{1/2}{n}^{-1}\right)+o\left(m^{-4}\right).
\end{eqnarray*}
To minimize the $MISE$ of $\widetilde{f}_{n,m,m/b}$, the order $(m_n)$ must equal to
\begin{eqnarray}
\left(\left[\frac{b^2}{\lambda_1(b)}\frac{8C_2}{C_1}\right]^{2/9}n^{2/9}\right),
\label{eq:36}
\end{eqnarray}
{and then the corresponding $MISE$}
\begin{eqnarray}
MISE\left[\widetilde{f}_{n,m,m/b}\right]=\left(b\lambda^4_1(b)\right)^{2/9}\frac{9\left(8C_1^8C_2\right)^{1/9}}{8}n^{-8/9}+o\left(n^{-8/9}\right).
\label{eq:37}
\end{eqnarray}
\begin{rem}
The equation (\ref{eq:37}) indicates that the choice $b=2$ is the best choice in terms of the $MISE$ for the estimator $\widetilde{f}_{n,m,m/b}$, since the factor $\left(b\lambda^4_1(b)\right)^{2/9}$ is increasing in $b=2,3, \ldots$.\\
This method of bias correction reduces the bias of Bernstein estimator from $O\left(m^{-1}\right)$ to $O\left(m^{-2}\right)$, but it loses the non-negativity. As an additive bias correction to the logarithm of estimator, \cite{Ter80} originally developed a multiplicative bias correction that enables keeping the non-negativity. This method was adopted by \cite{Hir10} for the beta kernel estimator introduced by \cite{Che99}. 
\end{rem}
\subsubsection{{Multiplicative bias-correction Bernstein density estimator $\widetilde{f}_{n,m,b,\varepsilon}$}}
This estimator was considered by~\cite{Kak14} by applying a multiplicative bias correction method to the Bernstein estimator
\begin{eqnarray}
\widetilde{f}_{n,m,b,\varepsilon}(x)=\left\{\widetilde{f}_{n,m}(x)\right\}^{b/(b-1)}\left\{\widetilde{f}_{n,m/b}(x)+\varepsilon\right\}^{-1/b-1}, \quad x\in[0,1],
\label{eq:38}
\end{eqnarray}
for some $\varepsilon=\varepsilon(m)>0$, converting to zero at a suitable rate. Under the Assumption $(A1)$, for $x\in[0,1]$ such as $f(x)>0$, with $m=O\left(n^\eta\right)$ and $\varepsilon\approx m^\tau$ where $\eta\in(0,1)$ and $\tau >2$  we have  
\begin{eqnarray*}
\mathbb{E}[\widetilde{f}_{n,m,b,\varepsilon}(x)]-f(x)&=&-\frac{b}{m^2}\left\{-\frac{\Delta^2_1(x)}{2f(x)}+\Delta_2(x)\right\}\\
&&+O\left(Var[\widetilde{f}_{n,m}(x)]+Var[\widetilde{f}_{n,m/b}(x)]\right)+o\left(m^{-2}\right),
\end{eqnarray*}
and
\begin{equation*}
Var[\widetilde{f}_{n,m,b,\varepsilon}(x)]=Var[\widetilde{f}_{n,m,m/b}(x)]+o\left(Var[\widetilde{f}_{n,m}(x)]+Var[\widetilde{f}_{n,m/b}(x)]+m^{-4}\right).
\end{equation*} 
Moreover, we have
\begin{eqnarray*}
MISE\left[\widetilde{f}_{n,m,b,\varepsilon}\right]=\lambda_1(b)C_1\frac{m^{1/2}}{n}+\frac{b^2C_5}{m^4}+o\left(m^{1/2}{n}\right)+o\left(m^{-4}\right).
\end{eqnarray*}
To minimize the $MISE$ of $\widetilde{f}_{n,m,b,\varepsilon}$, the order $\left(m_n\right)$ must equal to
\begin{eqnarray}
\left(\left[\frac{b^2}{\lambda_1(b)}\frac{8C_5}{C_1}\right]^{2/9}n^{2/9}\right),
\label{eq:39}
\end{eqnarray}
and the corresponding $MISE$ {is equal to}
\begin{eqnarray*}
MISE\left[\widetilde{f}_{n,m,b,\varepsilon}\right]=\left(b\lambda^4_1(b)\right)^{2/9}\frac{9\left(8C_1^8C_5\right)^{1/9}}{8}n^{-8/9}+o\left(n^{-8/9}\right).
\label{eq:40}
\end{eqnarray*}
\begin{rem}
Note that the estimator $\widetilde{f}_{n,m,b,\varepsilon}$ retains non-negativity, but it is not a genuine density. In fact, $\widetilde{f}_{n,m,b,\varepsilon}$ does not generally integrate to unity. To solve this problem, \cite{Kak14} proposed the normalized bias-corrected Bernstein estimator.
\end{rem}
\subsubsection{Normalized bias-corrected Bernstein estimator $\widetilde{f}^N_{n,m,b,\varepsilon}$}
The normalized bias-corrected Bernstein estimator is given by:
\begin{eqnarray}
\widetilde{f}^N_{n,m,b,\varepsilon}(x)=\frac{\widetilde{f}_{n,m,b,\varepsilon}(x)}{\displaystyle\int_{0}^1\widetilde{f}_{n,m,b,\varepsilon}(y)dy}, \quad x\in[0,1].
\label{eq:41}
\end{eqnarray}
Under the Assumption $(A1)$, for $x\in[0,1]$ such as $f(x)>0$, with $m=O\left(n^\eta\right)$ and $\varepsilon\approx m^\tau$ where $\eta\in(0,1)$ and $\tau >2$  we have  
\begin{eqnarray*}
\mathbb{E}[\widetilde{f}^N_{n,m,b,\varepsilon}(x)]-f(x)&=&-\frac{b}{m^2}\left\{-\frac{\Delta^2_1(x)}{2f(x)}+\Delta_2(x)+f(x)\displaystyle\int_0^1\frac{\Delta^2_1(y)}{2f(y)}dy\right\}\\
&&+O\left(Var[\widetilde{f}_{n,m}(x)]+Var[\widetilde{f}_{n,m/b}(x)]+n^{-1}m^{1/2}\right)+o\left(m^{-2}\right),
\end{eqnarray*}
and
\begin{equation*}
Var[\widetilde{f}^N_{n,m,b,\varepsilon}(x)]=Var[\widetilde{f}_{n,m,m/b}(x)]+o\left(Var[\widetilde{f}_{n,m}(x)]+Var[\widetilde{f}_{n,m/b}(x)]+n^{-1}m^{1/2}+m^{-4}\right),
\end{equation*} 
and then, we have
\begin{eqnarray*}
MISE\left[\widetilde{f}^N_{n,m,b,\varepsilon}\right]=\lambda_1(b)C_1\frac{m^{1/2}}{n}+\frac{b^2C_6}{m^4}+o\left(m^{1/2}{n}^{-1}\right)+o\left(m^{-4}\right).
\end{eqnarray*}
To minimize the $MISE$ of $\widetilde{f}^N_{n,m,b,\varepsilon}$, the order $\left(m_n\right)$ must equal to
\begin{eqnarray}
\left(\left[\frac{b^2}{\lambda_1(b)}\frac{8C_6}{C_1}\right]^{2/9}n^{2/9}\right),
\label{eq:42}
\end{eqnarray}
and the corresponding $MISE$ is equal to
\begin{eqnarray*}
MISE\left[\widetilde{f}^N_{n,m,b,\varepsilon}\right]=\left(b\lambda^4_1(b)\right)^{2/9}\frac{9\left(8C_1^8C_6\right)^{1/9}}{8}n^{-8/9}+o\left(n^{-8/9}\right).
\label{eq:43}
\end{eqnarray*}

\section{Applications}\label{section:app}
When using the Bernstein polynomials, we must consider a density on $[0,1]$. For this purpose, we need to make some suitable transformations in different cases (we list below) :
\begin{enumerate}
\item Suppose that $X$ is concentrated on a  finite support $[a,b]$. Then we work with the sample values $Y_1,\ldots, Y_n$ where $Y_i=\frac{X_i-a}{b-a}$. Denoting by $g_n(x)$ the estimated density function of $Y_1, \ldots, Y_n$, we compute the estimated density $f_n$ of $X$
\begin{eqnarray*}
f_n(x)=\frac{1}{b-a}g_n\left(\frac{x-a}{b-a}\right).
\end{eqnarray*}
\item For the densities functions concentrated on the interval $(-\infty,+\infty)$, we can use the transformed sample $Y_i=\frac{1}{2}+\frac{1}{\pi}\arctan(X_i)$ which transforms the range to the interval $(0,1)$. Hence 
\begin{eqnarray*}
f_n(x)=\frac{1}{\pi(1+x^2)}g_n\left(\frac{1}{2}+\frac{1}{\pi}\arctan(x)\right).
\end{eqnarray*}
\item For the support $[0,\infty)$, we can use the transformed sample $Y_i=\frac{X_i}{X_i+1}$ which  transforms the range to the interval $(0, 1)$. Hence 
\begin{eqnarray*}
f_n(x)=\frac{1}{(1+x)^2}g_n\left(\frac{x}{1+x}\right).
\end{eqnarray*}
\end{enumerate}

\paragraph{Computational cost} 
As mentioned in the introduction, the advantage of the proposed recursive estimators on their non-recursive version
is that their update, from a sample of size $n$ to one of size $n+1$, require less computations. This property can be generalized, one can check that it follows from~(\ref{eq:rec:density}) that for all $n_1\in \left[0,n-1\right]$,
\begin{eqnarray*}
f_{n}\left(x\right)&=&\prod_{j=n_1+1}^n\left(1-\gamma_{j}\right)f_{n-1}\left(x\right)+\sum_{k=n_1}^{n-1}\prod_{j=k+1}^n\left(1-\gamma_{j}\right)\gamma_{k}Z_k\left(x\right)+\gamma_nZ_n\left(x\right),\\
&=&\alpha_1f_{n-1}\left(x\right)+\sum_{k=n_1}^{n-1}\beta_{k}Z_k\left(x\right)+\gamma_nZ_n\left(x\right),
\end{eqnarray*}
where $\alpha_1=\prod_{j=n_1+1}^n\left(1-\gamma_{j}\right)$ and $\beta_k=\gamma_{k}\prod_{j=k+1}^n\left(1-\gamma_{j}\right)$.
{It is clear, that the proposed estimators can be viewed as a linear combination of two estimators, which improve considerably the computational cost.
Moreover, in order to give some comparative elements with Vitale's estimator defined in~\eqref{eq:26}, including computational coasts, we consider $500$ samples of size $500$ generated from the beta distribution $\mathcal{B}\left(3,5\right)$; moreover, we suppose that we receive an additional $500$ samples of size $500$ generated also from the beta distribution $\mathcal{B}\left(3,5\right)$. Performing the two methods, the running time using our proposed recursive estimator defined by algorithm~\eqref{eq:rec:density} with stepsize $\left(\gamma_n\right)=\left(n^{-1}\right)$ and the order $\left(m_n\right)$ according to the minimization of Least Squares Cross-Validation (LSCV) described below was roughly $25$ minutes on the author's workstation, while the running time using the non-recursive estimator~\eqref{eq:26} was roughly more than $4$ hours on the author's workstation. 
}
The aim of this paragraph is to compare the performance of Vitale's estimator $\widetilde{f}_n$ defined in (\ref{eq:26}), Leblanc's estimator $\widetilde{f}_{n,m,m/2}$ given in~(\ref{eq:32}), the generalized estimator $\widetilde{f}_{n,m,m/b}$, defined in~\eqref{eq:35}, the multiplicative bias corrected Bernstein estimator $\widetilde{f}_{n,m,b,\varepsilon}$ given in~(\ref{eq:38}) and the normalized estimator $\widetilde{f}^N_{n,m,b,\varepsilon}$ given in (\ref{eq:41}) with that of the proposed estimator (\ref{eq:1}).
\begin{itemize}
\item [(1)] When applying $f_n$, one needs to choose two quantities:
\begin{itemize}
\item [$\bullet$] The stepsize $(\g_n)$ {is chosen to be equal to $\left(\nu n^{-1}\right)$, with $\nu \in\left\{\frac{4}{5},\frac{8}{9},1\right\}$}.
\item [$\bullet$] The order $(m_n)$ is chosen to be equal to (\ref{eq:30}).
\end{itemize}
\item [(2)] When applying $\widetilde{f}_n$, one needs to choose the order $\left(m_n\right)$ to be equal to (\ref{eq:27}).
\item [(3)] When applying $\widetilde{f}_{n,m,m/2}$, one  needs to choose the order $\left(m_n\right)$ to be equal to (\ref{eq:33}).
\item[(3)] When applying $\widetilde{f}_{n,m,m/b}$, one  needs to choose the order $\left(m_n\right)$ to be equal to (\ref{eq:36}) and $b=3,4$.
\item [(4)] When applying $\widetilde{f}_{n,m,b,\varepsilon}$, one  needs to choose the order $\left(m_n\right)$ to be equal to (\ref{eq:39}), $b=2,3,4$  and $\varepsilon=0.00001$.
\item[(5)]  When applying $\widetilde{f}^N_{n,m,b,\varepsilon}$, one  needs to choose the order $\left(m_n\right)$ to be equal to (\ref{eq:42}), $b=2,3,4$  and $\varepsilon=0.00001$.
\end{itemize}
\subsection{Simulations}\label{subsection:simu}
We consider the following ten density functions :
\begin{description}
\item(a) the beta density $\mathcal{B}(3,5)$, $f(x)=\frac{x^2(1-x)^4}{B(3,5)}$,

\item(b) the beta density $\mathcal{B}(1,6)$, $f(x)=\frac{(1-x)^5}{B(1,6)}$,

\item(c) the beta density $\mathcal{B}(3,1)$, $f(x)=\frac{x^2}{B(3,1)}$,
 
\item(d) the beta mixture density $1/2\mathcal{B}(3,9)+1/2\mathcal{B}(9,3)$, $f(x)=0.5\frac{x^2(1-x)^8}{B(3,9)}+0.5\frac{x^8(1-x)^2}{B(9,3)}$,

\item(e) the beta mixture density $1/2\mathcal{B}(3,1)+1/2\mathcal{B}(10,10)$, $f(x)=0.5\frac{x^2}{B(3,1)}+0.5\frac{x^9(1-x)^9}{B(10,10)}$,  

\item (f)  the beta mixture density $1/2\mathcal{B}(1,6)+1/2\mathcal{B}(3,5)$,
$f(x)=0.5\frac{(1-x)^5}{B(1,6)}+0.5\frac{x^2(1-x)^4}{B(3,5)}$, 

\item (g) the beta mixture density $1/2\mathcal{B}(2,1)+1/2\mathcal{B}(1,4)$, $f(x)=0.5\frac{x}{B(2,1)}+0.5\frac{(1-x)^3}{B(1,4)}$,

\item (h) the truncated exponential density $\mathcal{E}_{[0,1]}(1/0.8)$, $f(x)=\frac{\exp(-x/0.8)}{0.8\left\{1-\exp(-1/0.8)\right\}}$,  

\item(i) the truncated normal density $\mathcal{N}_{[0,1]}(0,1)$, $f(x)=\frac{\exp(-x^2/2)}{\displaystyle\int_0^1\exp(-t^2/2)dt}$,

\item(j) the truncated normal mixture density $1/4\mathcal{N}_{[0,1]}(2,1)+3/4\mathcal{N}_{[0,1]}(-3,1)$,\\
$f(x)=0.25\frac{\exp(-(x-2)^2/2)}{\displaystyle\int_0^1\exp(-(t-2)^2/2)dt}+0.75\frac{\exp(-(x+3)^2/2)}{\displaystyle\int_0^1\exp(-(t+3)^2/2)dt}$.
\end{description}
For each density function and sample of size $n$, we approximate the
average integrated squared error ($ISE$) of the estimator using $N=500$ trials of sample size $n$; $\overline{ISE}=\frac{1}{N}\displaystyle\sum_{k=1}^{N}ISE\left[\hat{f}_k\right]$, where $\hat{f}_k$ is the estimator computed from the $k$th sample, and, $ISE[\hat{f}]=\displaystyle\int_{0}^1\left\{\hat{f}(x)-f(x)\right\}^2dx$.
\begin{table}
  \centering
  \begin{tabular}{c|c|c|c|c|c}
\hline Density function&n&Vitale's estimator&Recursive 1&Recursive 2&Recursive 3\\

\hline$(a)$&50&\textbf{0.085389}&0.088252&0.089350&0.092738\\
&200&0.028167&\textbf{0.025737}&0.026057&0.027045\\
&500&0.013533&\textbf{0.011398}&0.011540&0.011977\\

\hline$(b)$&50&0.145441&\textbf{0.129385}&0.130996&0.135963\\
&200&0.047977&\textbf{0.037733}&0.038202&0.039651\\
&500&0.023050&\textbf{0.016710}&0.016918&0.017560\\

\hline$(c)$&50&0.074634&\textbf{0.061163}&0.061925&0.064273\\
&200&0.024620&\textbf{0.017837}&0.018059&0.0187441\\
&500&0.011828&\textbf{0.007899}&0.007997&0.0083012\\

\hline$(d)$&50&\textbf{0.108664}&0.119492&0.120980&0.125567\\
&200&0.035845&\textbf{0.034847}&0.035281&0.036619\\
&500&0.017222&\textbf{0.015433}&0.015625&0.016217\\

\hline$(e)$&50&\textbf{0.124816}&0.157438&0.159398&0.165442\\
&200&\textbf{0.041174}&0.045914&0.046485&0.048248\\
&500&\textbf{0.019782}&0.020333&0.020587&0.021367\\

\hline$(f)$&50&\textbf{0.084119}&0.109237&0.110596&0.114790\\
&200&\textbf{0.027749}&0.031857&0.032253&0.033476\\
&500&\textbf{0.013332}&0.014108&0.014284&0.014825\\

\hline$(g)$&50&\textbf{0.063962}&0.065611&0.066428&0.068947\\
&200&0.021099&\textbf{0.019134}&0.019372&0.020107\\
&500&0.010137&\textbf{0.008474}&0.008579&0.008904\\

\hline$(h)$&50&\textbf{0.046147}&0.042890&0.043424&0.04507\\
&200&0.015222&\textbf{0.012508}&0.012664&0.013144\\
&500&0.007313&\textbf{0.005539}&0.005608&0.005821\\

\hline$(i)$&50&\textbf{0.031391}&0.036611&0.037066&0.038472\\
&200&\textbf{0.010355}&0.010676&0.010809&0.011219\\
&500&0.004975&\textbf{0.004728}&0.004787&0.004968\\

\hline$(j)$&50&0.071621&\textbf{0.065671}&0.066489&0.069010\\
&200&0.023626&\textbf{0.019152}&0.019390&0.020125\\
&500&0.011351&\textbf{0.008481}&0.008587&0.008913\\
\end{tabular}
\caption{The average integrated squared error ($ISE$) of \texttt{Vitale's estimator} $\widetilde{f}_n$ and the three recursive estimators; \texttt{recursive 1} correspond to the estimator~$f_n$ with the choice $(\gamma_n)=(n^{-1})$, \texttt{recursive 2} correspond to the estimator~$f_n$ with the choice $(\g_n)=\left(\left[1-\frac{a}{2}\right]n^{-1}\right)$ ($a=2/9$) and \texttt{recursive 3} correspond to the estimator~$f_n$ with the choice $(\g_n)=\left(\left[1-a\right]n^{-1}\right)$ ($a=1/5$).}
\label{Tab:1}
\end{table}

\begin{table}
  \centering
\small
\begin{tabular}{p{0.5cm}|c|c|c|c|c|c|c|c|c|c}
\multicolumn{3}{c|}{}&\multicolumn{2}{c|}{$\widetilde{f}_{n,m,b}$} &\multicolumn{3}{c|}{$\widetilde{f}_{n,m,b,\varepsilon}$, $\varepsilon=0.00001$}&\multicolumn{3}{c}{$\widetilde{f}^N_{n,m,b,\varepsilon}$, $\varepsilon=0.00001$}\\
\hline&n&$\widetilde{f}_{n,m,m/2}$&$b=3$&$b=4$&$b=2$&$b=3$&$b=4$&$b=2$&$b=3$&$b=4$\\
\hline$(a)$&50&{\bf 0.08504}&0.08687&0.08874&0.10597&0.10824&0.11057&0.09852&0.10064&0.10280\\
&200&{\bf 0.02480}&0.02533&0.02587&0.03090&0.03156&0.03224&0.02873&0.02935&0.02998\\
&500&{\bf 0.01098}&0.01122&0.01146&0.01368&0.01398&0.01428&0.01272&0.01299&0.01327\\

\hline$(b)$&50&{\bf 0.12469}&0.12736&0.13010&0.13284&0.13569&0.13861&0.14319&0.14626&0.14940\\
&200&{\bf 0.03636}&0.03714&0.03794&0.03874&0.03957&0.04042&0.04176&0.04265&0.04357\\
&500&{\bf 0.01610}&0.01644&0.01680&0.01715&0.01752&0.01790&0.01849&0.01889&0.01929\\

\hline$(c)$&50&{\bf 0.05894}&0.06020&0.06150&0.08278&0.08455&0.08637&0.08721&0.08908&0.09100\\
&200&{\bf 0.01719}&0.01755&0.01793&0.02414&0.02465&0.02518&0.02543&0.02598&0.02653\\
&500&{\bf 0.00761}&0.00777&0.00794&0.01069&0.01092&0.01115&0.01126&0.01150&0.01175\\

\hline$(d)$&50&{\bf 0.11515}&0.11762&0.12015&0.14391&0.14699&0.15015&0.14656&0.14971&0.15292\\
&200&{\bf 0.03358}&0.03430&0.03504&0.04196&0.04286&0.04379&0.04274&0.04366&0.04459\\
&500&{\bf 0.01487}&0.01519&0.01551&0.01858&0.01898&0.01939&0.01892&0.01933&0.01975\\

\hline$(e)$&50&{\bf 0.15172}&0.15498&0.15830&0.15311&0.15640&0.15976&0.15463&0.15795&0.16134\\
&200&{\bf 0.04424}&0.04519&0.04616&0.04465&0.04561&0.04659&0.04509&0.04606&0.04705\\ 
&500&{\bf 0.01959}&0.02001&0.02044&0.01977&0.02020&0.02063&0.01997&0.02040&0.02083\\

\hline$(f)$&50&0.10527&0.10753&0.10984&{\bf 0.10527}&0.11383&0.11627&0.11420&0.11665&0.11916\\
&200&{\bf 0.03070}&0.03135&0.03203&0.03250&0.03319&0.03391&0.03330&0.03402&0.03475\\
&500&{\bf 0.01359}&0.01388&0.01418&0.01439&0.01470&0.01501&0.01475&0.01506&0.01539\\

\hline$(g)$&50&0.06323&0.06458&0.06597&0.06298&0.06433&0.06571&{\bf 0.06252}&0.06386&0.06524\\
&200&0.01844&0.01883&0.01924&0.01836&0.01876&0.01916&{\bf 0.01823}&0.01862&0.01902\\
&500&0.00816&0.00834&0.00852&0.00813&0.00830&0.00848&{\bf 0.00807}&0.00824&0.00842\\

\hline$(h)$&50&0.04133&0.04222&0.04312&0.03872&0.03955&0.04040&{\bf 0.03566}&0.03643&0.03721\\ 
&200&0.01205&0.01231&0.01257&0.01129&0.01153&0.01178&{\bf 0.01040}&0.01062&0.01085\\ 
&500&0.00533&0.00545&0.00557&0.00500&0.00510&0.00521&{\bf 0.00460}&0.00470&0.00480\\

\hline$(i)$&50&{\bf 0.03528}&0.03603&0.03681&0.03594&0.03671&0.03750&0.03589&0.03666&0.03745\\
&200&{\bf 0.01028}&0.01051&0.01073&0.01048&0.01070&0.01093&0.01046&0.01069&0.01092\\
&500&{\bf 0.00455}&0.00465&0.00475&0.00464&0.00474&0.00484&0.00463&0.00473&0.00483\\

\hline$(j)$&50&0.06328&0.06464&0.06603&0.06439&0.06577&0.06718&{\bf 0.06236}&0.06370&0.06507\\ 
&200&0.01845&0.01885&0.01925&0.01877&0.01918&0.01959&{\bf 0.01818}&0.01857&0.01897\\
&500&0.00817&0.00817&0.00852&0.00831&0.00849&0.00867&{\bf 0.00805}&0.00822&0.00840\\
\end{tabular}
\caption{The average integrated squared error ($ISE$) of Leblanc estimator's $\widetilde{f}_{n,m,m/2}$ and the three estimators introduced by Kakizawa: $\widetilde{f}_{n,m,m/b}$ with $b=3,4$, $\widetilde{f}_{n,m,b,\varepsilon}$ and $\widetilde{f}^N_{n,m,b,\varepsilon}$ with $b=2,3,4$ and $\varepsilon=0.00001$.}
\label{Tab:2}
\end{table}
From Table \ref{Tab:1} and \ref{Tab:2} we conclude that : 
\begin{itemize}
\item In all the considered cases, the average $ISE$ of our density estimator (\ref{eq:1}) is smaller than that of Vitale's estimator defined in (\ref{eq:26}), except the cases $(e)$ and $(f)$ of the Beta mixture and the cases $(a)$, $(d)$, $(g)$ and $(h)$ for the small sample size $n=50$ and in the case $(i)$ for the size $n=50$ and $n=200$.
\item In all the cases, the average $ISE$ of our recursive density estimator (\ref{eq:1}) is slightly larger then that of Leblanc's estimator $\widetilde{f}_{n,m,m/2}$ given in (\ref{eq:32}).

\item In all the cases, the average $ISE$ of our density estimator (\ref{eq:1}) is smaller than that of the generalized estimator $\widetilde{f}_{n,m,m/b}$, with $b=4$ (see (\ref{eq:35})).
 
\item In all the cases, the average $ISE$ of our density estimator (\ref{eq:1}) is smaller than that of the multiplicative bias corrected Bernstein estimator $\widetilde{f}_{n,m,b,\varepsilon}$ given in~(\ref{eq:38}) and the normalized estimator $\widetilde{f}^N_{n,m,b,\varepsilon}$ given in~(\ref{eq:41}), except the cases $(g)$ and $(h)$.

\item The average $ISE$ of the generalized estimator $\widetilde{f}_{n,m,m/b}$, (\ref{eq:35}) increase when $b$ increase, {then} the optimal choice is $b=2$ which corresponds to Leblanc's estimator $\widetilde{f}_{n,m,m/2}$ given in~(\ref{eq:32}).

\item The average $ISE$ of the multiplicative bias corrected Bernstein estimator $\widetilde{f}_{n,m,b,\varepsilon}$ given in~(\ref{eq:38}) and the normalized estimator $\widetilde{f}^N_{n,m,b,\varepsilon}$ defined in~(\ref{eq:41}) increase when $b$ increase

\item The average $ISE$ of the multiplicative bias corrected Bernstein estimator $\widetilde{f}_{n,m,b,\varepsilon}$ given in~(\ref{eq:38}) and the normalized estimator $\widetilde{f}^N_{n,m,b,\varepsilon}$ given in~(\ref{eq:41}) are smaller than that of the estimator $\widetilde{f}_{n,m,m/b}$ defined in~(\ref{eq:35}), in the cases $(g)$ and $(h)$ and are larger in the other cases.

\item The average $ISE$ of the multiplicative bias corrected Bernstein estimator $\widetilde{f}_{n,m,b,\varepsilon}$ (\ref{eq:38}) is larger than that of the normalized estimator $\widetilde{f}^N_{n,m,b,\varepsilon}$ given in~(\ref{eq:41}) in the cases $(a)$, $(g)$, $(h)$, $(i)$ and $(j)$ and is smaller in the other cases.

\item The average $ISE$ decreases as the sample size increases.
\end{itemize}
\begin{figure}[!ht]
\begin{center}
\includegraphics[width=0.6\textwidth,angle=270,clip=true,trim=40 0 0 
0]{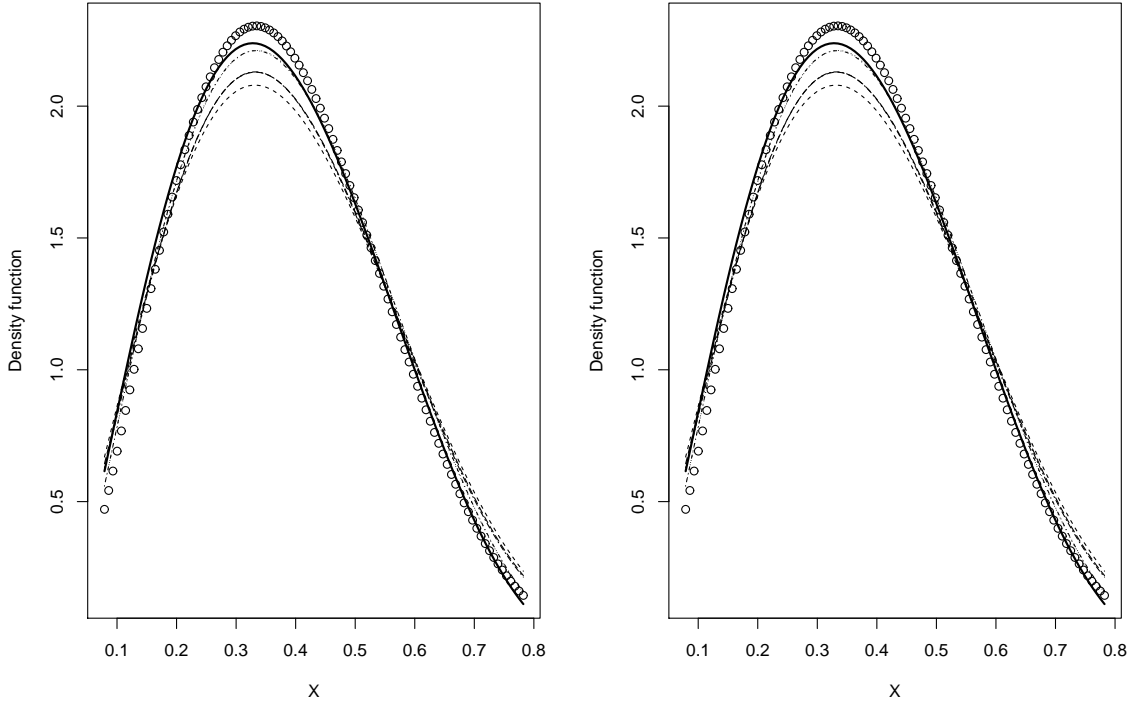}
\end{center}
\label{Fig:1}
\caption{Qualitative comparison between the proposed density estimator $f_n$ given in~(\ref{eq:1}) with stepsize $(\gamma_n)=(n^{-1})$ (solid line), the Vitale's estimator $\widetilde{f}_n$ defined in~(\ref{eq:26}) (dashed line), Leblanc's estimator $\widetilde{f}_{n,m,m/2}$ defined in~(\ref{eq:32}) (dotted line), the generalized estimator $\widetilde{f}_{n,m,m/4}$ defined in~(\ref{eq:35}) (dot-dashed line), Kakizwa's estimator $\widetilde{f}_{n,m,b,\varepsilon}$ defined in~(\ref{eq:38}) (long-dashed) and the normalized estimator $\widetilde{f}^N_{n,m,b,\varepsilon}$ defined in~(\ref{eq:41}) (two-dashed line) with $b=2$ and $\varepsilon=0.00001$ for $500$ samples respectively of size $50$ (left panel) and of size $250$ (right panel) of the beta density $\mathcal{B}\left(3,5\right)$.}
\end{figure}
From figure \ref{Fig:1}, we can observe that :
\begin{itemize}
\item The proposed recursive estimator $f_n$ given in~(\ref{eq:1}) using the stepsize $\left(\gamma_n\right)=\left(n^{-1}\right)$ is closer to the true  density function compared to Vitale's estimator $\widetilde{f}_n$ given in (\ref{eq:26}), the generalized estimator $\widetilde{f}_{n,m,m/4}$ given in~(\ref{eq:35}), Kakizwa's estimator $\widetilde{f}_{n,m,b,\varepsilon}$ defined in~(\ref{eq:38}) and the normalized estimator $\widetilde{f}^N_{n,m,b,\varepsilon}$ defined in~(\ref{eq:41}) with $b=2$ and $\varepsilon=0.00001$.
\item Within the interval $[0,1]$, our density estimator (\ref{eq:1}) using the stepsize $\left(\gamma_n\right)=\left(n^{-1}\right)$ is closer to the true  density function compared to Leblanc's estimator $\widetilde{f}_{n,m,m/2}$ defined in~(\ref{eq:32}).
\item When the sample size increases, we get closer estimation of the true density function.
\end{itemize}
\subsection{Real dataset}\label{subsection:realdata}
In any practical situation, to estimate an unknown density function $f$, it is essential to specify the order $m$  to be used for the estimator. {One way is to use least squares cross-validation} ($LSCV$) method to obtain a data-driven choice of $m$.
\subsubsection{Order selection method}
First, we recall that the $LSCV$ method is based on minimizing the
integrated square error between the estimated density function $\hat{f}$ and
the true density function $f$
\begin{eqnarray*}
\displaystyle\int_0^1 \left(\hat{f}(x)-f(x)\right)^2dx=\displaystyle\int_0^1\hat{f}^2(x)dx-2\displaystyle\int_0^1\hat{f}(x)f(x)dx+\displaystyle\int_0^1f^2(x)dx.
\end{eqnarray*}
From this, \cite{Sil86} derived the score function 
\begin{eqnarray}
LSCV_{\hat{f}}(m)=\displaystyle\int_0^1\hat{f}^2(x)dx-\frac{2}{n}\displaystyle\sum_{i=1}^{n}\hat{f}_{-i}(X_i),
\label{eq:45}
\end{eqnarray}
where $\hat{f}_{-i}$ is the density
estimate without the data point $X_i$. The smoothing parameter is chosen by minimizing $LSCV(m)$ ($m=\arg\min_{m}LSCV(m)$). \\
For our proposed recursive estimator $f_n$ given in (\ref{eq:1}), we make the following choice of $(\g_n)=(n^{-1})$, and then we choice the order $m$ in order to minimize the following criterion :
\begin{eqnarray*}
LSCV_{f_n}(m)=\displaystyle\int_0^1f^2_n(x)dx-\frac{2}{n}\displaystyle\sum_{i=1}^{n}f_{n,-i}(X_i).
\end{eqnarray*}
We define integer sequences $p_i=[m_iX_i]$ and $q_i=[m_iX_i/2]$, {we then have} $X_i\in\left(\frac{p_i}{m_i}, \frac{p_{i}+1}{m_i}\right]$ and $X_i\in\left(\frac{2q_i}{m_i}, \frac{2(q_i+1)}{m_i}\right]$. Then we obtain 
\begin{eqnarray*}
f_{n,-i}(x)=\frac{1}{n-1}\left[nf_n(x)-\displaystyle\sum_{i=1}^{n}\left\{2m_ib_{p_i}(m_i-1,x)-\frac{m_i}{2}b_{q_i}\left(\frac{m_i}{2}-1,x\right)\right\}\right],
\end{eqnarray*}
{then} we conclude that
\begin{eqnarray*}
LSCV_{f_n}(m)=\displaystyle\int_0^1f^2_n(x)dx-\frac{2}{n-1}\left[\displaystyle\sum_{i=1}^{n}f_{n}(X_i)-\frac{1}{n}\displaystyle\sum_{i=1}^{n}\left\{2m_ib_{p_i}(m_i-1,X_i)-\frac{m_i}{2}b_{q_i}\left(\frac{m_i}{2}-1,X_i\right)\right\}\right].
\end{eqnarray*} 
Note that, the $LSCV$ function for Vitale's  estimator (\ref{eq:26}), is written as 
\begin{eqnarray*}
LSCV_{\widetilde{f}_n}(m)=\displaystyle\int_0^1\left\{\widetilde{f}_n(x)\right\}^2dx-\frac{2}{n}\displaystyle\sum_{i=1}^{n}\widetilde{f}_{n,-i}(X_i).
\end{eqnarray*}
We define integer sequences $k_i=[mX_i]$, {we then have} $X_i\in\left(\frac{k_i}{m}, \frac{k_i+1}{m}\right]$, {consequently} 
\begin{eqnarray*}
\widetilde{f}_{n,-i}(x)&=&\frac{1}{n-1}\left[\widetilde{f}_n(x)-mb_{k_i}(m-1,x)\right],
\end{eqnarray*}
{then} we conclude 
\begin{eqnarray*}
LSCV_{\widetilde{f}_n}(m)=\displaystyle\int_0^1\left\{\widetilde{f}_n(x)\right\}^2dx-\frac{2}{n-1}\left[\displaystyle\sum_{i=1}^n\widetilde{f}_n(X_i)-\frac{m}{n}\displaystyle\sum_{i=1}^nb_{k_i}(m-1,X_i)\right].
\end{eqnarray*} 
The $LSCV$ function for the estimator $\widetilde{f}_{n,m,m/b}$ defined in (\ref{eq:35}) with $b=2,3,4,...$, is written as 
\begin{eqnarray*}
LSCV_{\widetilde{f}_{n,m,m/b}}(m)=\displaystyle\int_0^1\left\{\widetilde{f}_{n,m,m/b}(x)\right\}^2dx-\frac{2}{n}\displaystyle\sum_{i=1}^{n}\widetilde{f}_{n,m,m/b,-i}(X_i).
\end{eqnarray*}
We define integer sequences $k_i=[mX_i]$ and $r_i=[mX_i/b]$, {we then have} $X_i\in\left(\frac{k_i}{m}, \frac{k_i+1}{m}\right]$ and $X_i\in\left(\frac{br_i}{m}, \frac{b(r_{i}+1)}{m}\right]$. Then we obtain 
\begin{eqnarray*}
LSCV_{\widetilde{f}_{n,m,m/b}}(m)&=&\displaystyle\int_0^1\left\{\widetilde{f}_{n,m,m/b}(x)\right\}^2dx-\frac{2}{n-1}\Bigg[\displaystyle\sum_{i=1}^{n}\widetilde{f}_{n,m,m/b}(X_i)\\
&-&\frac{1}{n}\displaystyle\sum_{i=1}^{n}\left\{\frac{b}{b-1}mb_{k_i}(m-1,X_i)-\frac{1}{b-1}\frac{m}{b}b_{r_i}\left(\frac{m}{b}-1,X_i\right)\right\}\Bigg].
\end{eqnarray*} 
Using the Kakizawa's estimators $\widetilde{f}_{n,m,b,\varepsilon}$ defined in~(\ref{eq:38}) and $\widetilde{f}^N_{n,m,b,\varepsilon}$ defined in~(\ref{eq:41}), the $LSCV$ function is written as in (\ref{eq:45}).\\

\subsubsection{Old Faithful data}\label{data:1}
\begin{figure}[!ht]
\begin{center}
\includegraphics[width=0.65\textwidth,angle=270,clip=true,trim=40 0 0 0]{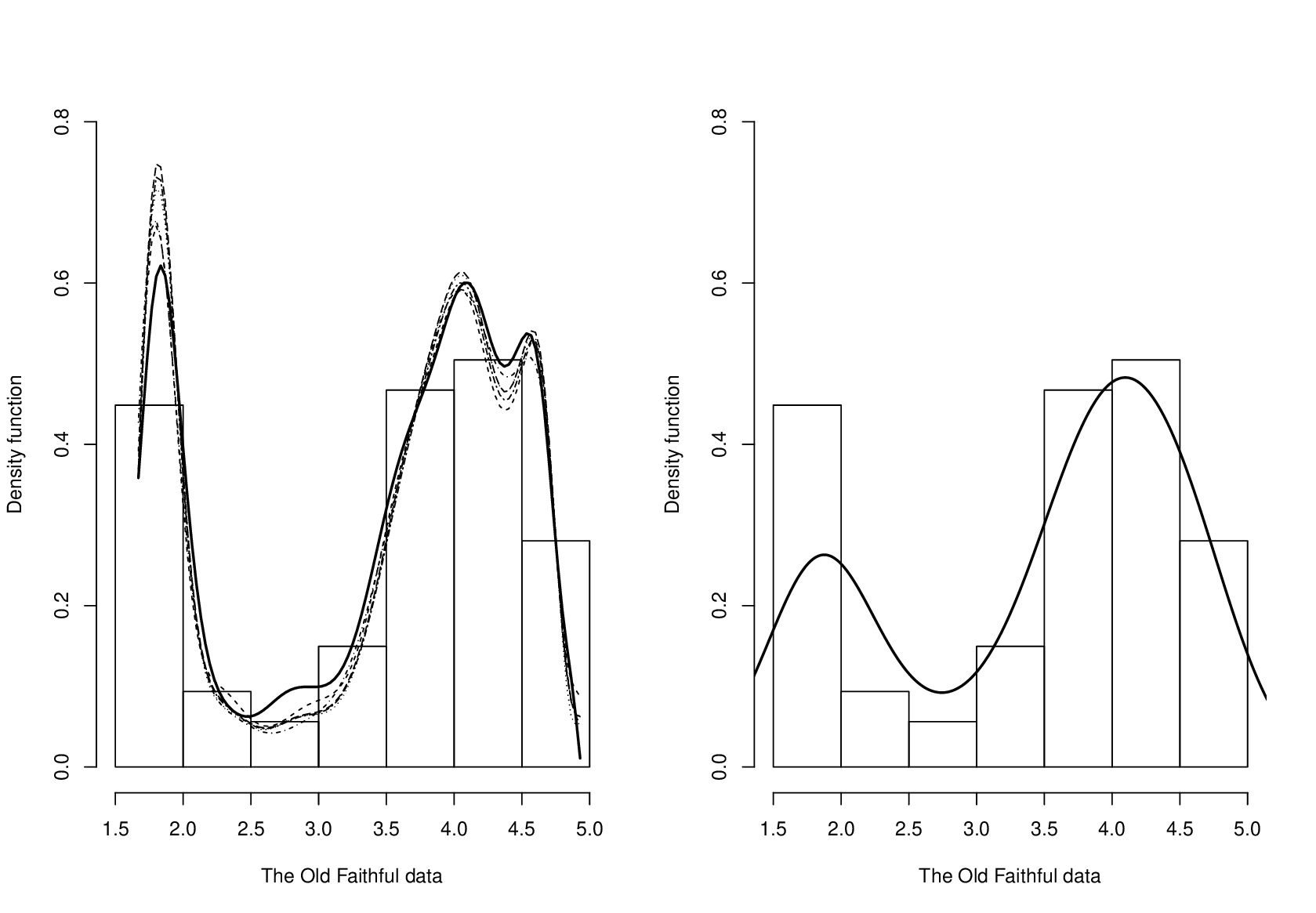}
\end{center}
\label{Fig:3}
\caption{Density estimates for the Old Faithful data : recursive estimator $f_n$ defined in~(\ref{eq:1}) (solid line), Vitale's estimator $\widetilde{f}_n$ defined in~(\ref{eq:26}) (dashed line), Leblanc's estimator defined in~$\widetilde{f}_{n,m,m/2}$ defined in~(\ref{eq:32}) (dotted line), the generalized estimator $\widetilde{f}_{n,m,m/4}$ defined in~(\ref{eq:35}) (dot-dashed line), Kakizwa's estimator $\widetilde{f}_{n,m,b,\varepsilon}$ defined in~(\ref{eq:38}) (long-dashed) and the normalized estimator $\widetilde{f}^N_{n,m,b,\varepsilon}$ defined in~(\ref{eq:41}) (two-dashed line) with $b=2$ and $\varepsilon=0.00001$  (left panel) and Gaussian kernel density estimate using the $LSCV$-based bandwidth $h=0.3677$ (right panel)} 
\end{figure}

In this subsection, we consider the well known Old Faithful data given in Tab.2.2 of \cite{Sil86}. These data consist of the eruption lengths (in minutes) of $107$ eruptions of the Old Faithful geyser in Yellowstone National Park, U.S.A. The data are such that $\min_i(x_i)=1.67$ and $\max_i(x_i)=4.93$, {then} it is convenient to assume that the density of eruption times is defined on the interval $[1.5,5]$ and transform the data into the unit interval.\\
The $LSCV$ procedure was performed and resulted in $m=104$ for Vitale's estimator $\widetilde{f}_n$ defined in~(\ref{eq:26}), $(m_n)=(n^{0.987})$ for our proposed estimator $f_n$ defined in~(\ref{eq:1}), $m=66$ for Leblanc's estimator $\widetilde{f}_{n,m,m/2}$ defined in~(\ref{eq:32}), $m=52$ for the estimator $\widetilde{f}_{n,m,m/4}$ defined in~(\ref{eq:35}), $m=66$ for the multiplicative bias corrected Bernstein estimator $\widetilde{f}_{n,m,2,0.00001}$ defined in~(\ref{eq:38}) and $m=66$ for the normalized estimator $\widetilde{f}^N_{n,m,2,0.00001}$ defined in~(\ref{eq:41}). These estimators are shown in Figure \ref{Fig:3} along with an histogram of the data and a Gaussian kernel density estimate using the $LSCV$-based bandwidth $h=0.3677$. All the estimators are smooth and seem to capture the pattern highlighted by the histogram. We observe that our recursive estimator outperformed the others estimators near $x=1.5$.
\subsubsection{Tuna data}\label{data:2}
\begin{figure}[!ht]
\begin{center}
\includegraphics[width=0.65\textwidth,angle=270,clip=true,trim=40 0 0 0]{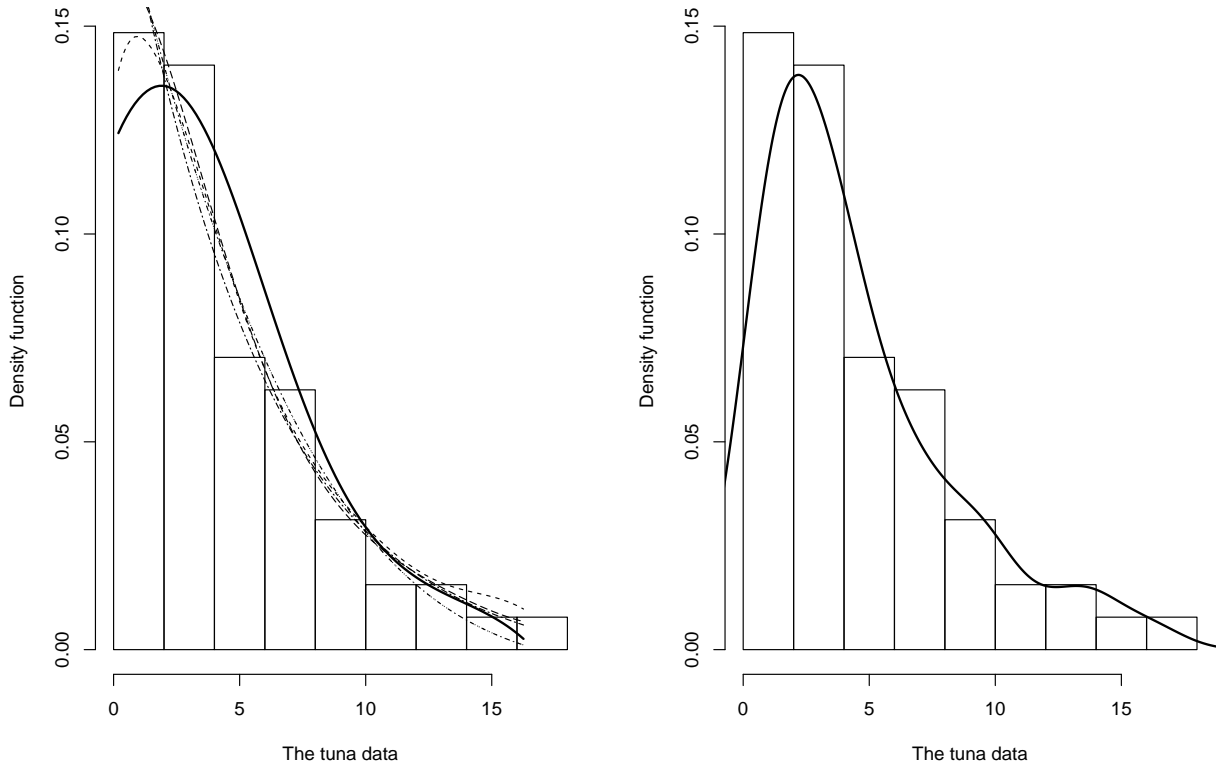}
\end{center}
\label{Fig:4}
\caption{Density estimates for the tuna data : recursive estimator $f_n$ defined in (\ref{eq:1}) (solid line), Vitale's estimator $\widetilde{f}_n$ defined in (\ref{eq:26}) (dashed line), Leblanc's estimator $\widetilde{f}_{n,m,m/2}$ defined in (\ref{eq:32}) (dotted line), the generalized estimator $\widetilde{f}_{n,m,m/4}$ defined in (\ref{eq:35}) (dot-dashed line), Kakizwa's estimator $\widetilde{f}_{n,m,b,\varepsilon}$ defined in (\ref{eq:38}) (long-dashed) and the normalized estimator $\widetilde{f}^N_{n,m,b,\varepsilon}$ defined in (\ref{eq:41}) (two-dashed line) with $b=2$ and $\varepsilon=0.00001$ (left panel) and Gaussian kernel density estimate using the $LSCV$-based bandwidth $h=1.291$ (right panel)} 
\end{figure}

Our last example concerne the tuna data given in \cite{Che96}. The data come from an aerial line transect survey of Southern Bluefin Tuna in the Great Australian Bight. An aircraft with two spotters on board flies randomly allocated line transects. The data are the perpendicular sighting distances (in miles) of $64$ detected tuna schools to the transect lines. The survey was conducted in summer when tuna tend to stay on the surface. We analyzed the transformed data divided by $w=18$ (the data are such that $\min_i(x_i)=0.19$ and $\min_i(x_i)=16.26$).\\
The $LSCV$ procedure was performed and resulted in $m=14$ for Vitale's estimator $\widetilde{f}_n$ defined in~(\ref{eq:26}), $(m_n)=(n^{0.633})$ for our proposed estimator $f_n$ defined in~(\ref{eq:1}), $m=4$ for Leblanc's estimator $\widetilde{f}_{n,m,m/2}$ defined in~(\ref{eq:32}), $m=4$ for the estimator $\widetilde{f}_{n,m,m/4}$ defined in~(\ref{eq:35}), $m=8$ for the multiplicative bias corrected Bernstein estimator $\widetilde{f}_{n,m,2,0.00001}$ defined in~(\ref{eq:38}) and $m=4$ for the normalized estimator $\widetilde{f}^N_{n,m,2,0.00001}$ defined in~(\ref{eq:41}). These estimators are shown in Figure \ref{Fig:4} along with an histogram of the data and a Gaussian kernel density estimate using the $LSCV$-based bandwidth $h= 1.291$. All the estimators are smooth and seem to capture the pattern highlighted by the histogram. We can observe that our proposed recursive estimator outperformed the other estimators near the boundaries. 
\section{Conclusion}\label{section:conclusion}
In this paper, we propose a recursive estimator of a density function based on a stochastic algorithm derived from Robbins-Monro's scheme and using Bernstein polynomials. We first study its asymptotic properties. We show that our proposed estimator of density function have a good boundary properties. Moreover, the bias rate of the proposed estimator is of $m^{-2}$, which is better than the Vitale's estimator with a bias rate of $m^{-1}$. For almost all the cases, the average $ISE$ of the proposed estimator (\ref{eq:1}) with a stepsize $\left(\g_n\right)=\left(n^{-1}\right)$ and the corresponding order $\left(m_n\right)$ is smaller than that of Vitale's estimator and than that of the multiplicative bias-corrected estimator defined by Kakizawa. Furthermore, our proposed recursive density estimator has a slightly larger average $ISE$ compared to Leblanc's estimator. In addition, a major advantage of our proposal is that its update, when new sample points are available, require less computational cost than the non recursive estimators. Our proposed estimator always integrates to unity, but is not necessarily non negative. However, we found that truncation and renormalisation may solve this issue. Finally, through simple real-life examples (Old Faithful data and tuna data) and a simulation study, we demonstrated  how the recursive Bernstein polynomial density estimators can lead to very satisfactory estimates. 
In conclusion, using the proposed recursive estimator $f_n$, we can obtain better results compared to those given by Vitale's estimator, Leblanc's estimator and the multiplicative bias-corrected estimator defined by Kakizawa especially near the boundaries.
\appendix
\section{Proofs}\label{section:proof}
In this section, we present proofs for the results presented in Section \ref{section:main:result}. We need the following technical lemma, which is proved in \cite{Mok09}.
\begin{lem}\label{lem:1}
Let $(v_{n})\in \mathcal{GS}(v^{*}),(\gamma_{n})\in \mathcal{GS}(-\alpha)$, and $l>0$ such that $l-v^{*}\xi>0$. We have\\
\begin{equation*}
\lim_{n\rightarrow\infty}v_{n}\Pi^{l}_{n}\sum^{n}_{k=1}\Pi^{-l}_{k}\frac{\gamma_{k}}{v_{k}}=\frac{1}{l-v^{*}\xi}.
\end{equation*}
Moreover, for all positive sequence $(\alpha_{n})$ such that $\lim_{n\rightarrow\infty}\alpha_{n}=0$, and all $\delta\in\mathbb{R}$, we have
\begin{equation*}
\lim_{n\rightarrow\infty}v_{n}\Pi^{l}_{n}\left[\sum^{n}_{k=1}\Pi^{-l}_{k}\frac{\gamma_{k}}{v_{k}}\alpha_{k}+\delta\right]=0.
\end{equation*}
\end{lem}
{Lemma~\ref{lem:1} is widely applied throughout the proofs. Let us underline that it is its application, which requires Assumption $\left(A4\right)$ on the limit of $\left(n\gamma_n\right)$ as $n$ goes to infinity. Let us mention that,} $\lim_{n\rightarrow\infty}(n\gamma_n)<\infty $ only if $\alpha =1$, the condition $\lim_{n\rightarrow\infty}(n\gamma_{n})\in(\min\left(a,(2\alpha+a)/4\right),\infty]$ in $\left(A4\right)$, which appears throughout our proofs, is equivalent to the condition $\lim_{n\rightarrow\infty}(n\gamma_{n})\in(\min\left(a,(2+a)/4\right),\infty]$. Similarly, since $\xi\neq 0$ only if $ \alpha=1$, we can consider $\alpha=1$ in all the results given in this paper.\\
\subsection{Proof of Proposition~\ref{prop:1} and Proposition~\ref{prop:2}}\label{sub:6.1}
In view of (\ref{eq:1}), we have
\begin{eqnarray}
\mathbb{E}\left[f_{n}(x)\right]-f(x)=\Pi_{n}\sum^{n}_{k=1}\Pi^{-1}_{k}\gamma_{k}\left(\mathbb{E}\left[Z_{k}(x)\right]-f(x)\right).
\label{eq:44}
\end{eqnarray}
Leblanc showed that, for $x\in[0,1]$, we have (see Theorem 6 in \cite{Leb10}) 
\begin{eqnarray*}
\mathbb{E}\left[Z_{n}(x)\right]-f(x)=-2\frac{\Delta_2(x)}{m_n^2}+o(m_n^2).
\end{eqnarray*}
Substituting this result into (\ref{eq:44}) leads to
\begin{eqnarray*}
\mathbb{E}\left[f_{n}(x)\right]-f(x)=\Pi_{n}\sum^{n}_{k=1}\Pi^{-1}_{k}\gamma_{k}\left[\frac{-2\Delta_2(x)}{m_k^2}(1+o(1)\right], \quad x\in [0,1].
\end{eqnarray*}
For $x\in(0,1)$, we obtain
\begin{itemize}
\item In the case ${a\in \left(0,\frac{2}{9}\alpha\right]}$, we have $\lim_{n\rightarrow\infty}(n\gamma_{n})>2a$. {Then the application of} lemma \ref{lem:1} gives
(\ref{eq:2}).
\item  In the case  ${a\in \left(\frac{2}{9}\alpha,1\right)}$, we have $m^{-2}_{n}=o\left(\sqrt{\gamma_{n}m^{1/2}_{n}}\right)$, since $2a>(2\alpha-a)/4$. {Then the application of} lemma \ref{lem:1} gives (\ref{eq:3}).
\end{itemize}
For $x=0, 1$, we obtain 
\begin{itemize}
\item In the case ${a\in \left(0,\frac{\alpha}{5}\right]}$, we have $\lim_{n\rightarrow\infty}(n\gamma_{n})>2a$. {Then the application of} lemma \ref{lem:1} gives (\ref{eq:2'}).
\item  In the case ${a\in \left(\frac{\alpha}{5},1\right)}$, we have $m^{-2}_{n}=o\left(\sqrt{\gamma_{n}m_{n}}\right)$ and  $2a>(\alpha-a)/2$. {Then the application of} lemma \ref{lem:1} gives (\ref{eq:3'}).
\end{itemize}
On the other hand, we have
\begin{eqnarray}
Var[f_{n}(x)]=\Pi^{2}_{n}\sum^{n}_{k=1}\Pi^{-2}_{k}\gamma^{2}_{k}Var[Z_{k}(x)].
\nonumber 
\end{eqnarray}
Leblanc derived that (see Theorem 6 in \cite{Leb10})
\begin{eqnarray*}
Var\left[Z_n(x)\right]=\begin{cases}
\left(\frac{1}{\sqrt{2}}+4\left(1-\sqrt{\frac{2}{3}}\right)\right)m_n^{1/2}f(x)\psi(x)+o_x(m_n^{1/2}),& \quad\text{for }x\in(0,1),\\
\frac{5}{2}m_nf(x)+o(m_n),&\quad\text{for }x=0, 1.
\end{cases}
\end{eqnarray*}
It follows that 
\begin{eqnarray*}
Var\left[f_n(x)\right]=\begin{cases}
C_3f(x)\psi(x)\Pi^{2}_{n}\sum^{n}_{k=1}\Pi^{-2}_{k}\gamma^{2}_{k}m_k^{1/2}(1+o_x(1)),& \quad\text{for }x\in(0,1),\\
\frac{5}{2}f(x)\Pi^{2}_{n}\sum^{n}_{k=1}\Pi^{-2}_{k}\gamma^{2}_{k}m_k(1++o(1)),&\quad\text{for }x=0, 1.
\end{cases}
\end{eqnarray*}
{Now, using} this result, for $x\in(0,1)$, we obtain 
\begin{itemize}
\item In the case {$a\in \left[\frac{2}{9}\alpha,1\right)$}, we have $\lim_{n\rightarrow\infty}(n\gamma_{n})>(2\alpha-a)/4$, the application of lemma \ref{lem:1} then gives
\begin{eqnarray*}
Var[f_{n}(x)]=C_3\gamma_{n}m_n^{1/2}\frac{2}{4-(2\alpha-a)\xi}f(x)\psi(x)+o\left(\gamma_{n}m_n^{1/2}\right), 
\end{eqnarray*}
and (\ref{eq:4}) follows.\\
\item  In the case  {$a\in \left(0,\frac{2}{9}\alpha\right)$}, we have $2a<(2\alpha-a)/4$ and $\gamma_{n}m_n^{1/2}=o\left(m^{-4}_{n}\right)$, {the application of lemma}~\ref{lem:1} ensures that
\begin{eqnarray*}
Var[f_{n}(x)]&=&\Pi_{n}\sum^{n}_{k=1}\Pi^{-1}_{k}\gamma_{k}o\left(m^{-4}_{k}\right),\nonumber\\
&=&o\left(m^{-4}_{n}\right),
\end{eqnarray*}
which gives (\ref{eq:5}).
\end{itemize}
Similarly, for $x=0,1$, we have 
\begin{itemize}
\item In the case {$a\in \left[\frac{\alpha}{5},1\right)$}, we have $\lim_{n\rightarrow\infty}(n\gamma_{n})>(\alpha-a)/2$, {the application of lemma}~\ref{lem:1} gives
(\ref{eq:4'}).
\item  In the case  {$a\in \left(0,\frac{\alpha}{5}\alpha\right)$}, we have $2a<(\alpha-a)/2$ and $\gamma_{n}m_n=o\left(m^{-4}_{n}\right)$, {the application of lemma}~\ref{lem:1} gives (\ref{eq:5'}).
\end{itemize}

\subsection{Proof of Proposition~\ref{prop:3}} \label{sub:6.2}
{Proposition \ref{prop:1} ensures that}
\begin{itemize}
\item [$\bullet$] In the case {$a\in \left(0,\frac{2}{9}\alpha\right)$}, we have
\begin{eqnarray}
\int^{1}_{0}\left\{Bias\left[f_n(x)\right]\right\}^{2}dx&=&\int^{1}_{0}\left[-2m^{-2}_{n}\frac{1}{1-2a\xi}\Delta_2(x)+o\left(m^{-2}_{n}\right)\right]^{2}dx\nonumber\\
&=&4C_2m^{-4}_{n}\frac{1}{(1-2a\xi)^{2}}+o\left(m^{-4}_{n}\right).
\label{eq:18}
\end{eqnarray}
\item [$\bullet$] In the case  {$a\in \left(\frac{2}{9}\alpha,1\right)$}, we have
\begin{eqnarray}
\int^{1}_{0}\left\{Bias\left[f_n(x)\right]\right\}^{2}dx=o\left(\gamma_{n}m^{1/2}_{n}\right).
\label{eq:19}
\end{eqnarray}
\end{itemize}
{Moreover}, we note that
\begin{eqnarray*}
\int^{1}_{0}Var[f_{n}(x)]dx=\Pi^{2}_{n}\sum^{n}_{k=1}\Pi^{-2}_{k}\gamma^{2}_{k}\int^{1}_{0}Var[Z_{k}(x)]dx.
\end{eqnarray*}
Since we have (see the proof of Theorem 7 in \cite{Leb10})
\begin{equation*}
\displaystyle\int_0^1 Var\left[Z_k(x)\right]dx=C_1C_3m_k^{1/2}+o(m_k^{1/2}),
\end{equation*}
{then we get}
\begin{eqnarray*}
\displaystyle\int_0^1Var\left[f_n(x)\right]=C_1C_3\Pi_n^2\displaystyle\sum_{k=1}^{n}\Pi_k^{-2}\gamma_k^2m_k^{1/2}(1+o(1)).
\end{eqnarray*}
\begin{itemize}
\item [$\bullet$] In the case {$a\in \left[\frac{2}{9}\alpha,1\right)$}, we have $\lim_{n\to+\infty}\left(n\gamma_n\right)>(2\alpha-a)/4$. {Then the application of} lemma \ref{lem:1} gives 
\begin{eqnarray}
\int^{1}_{0}Var\left[f_n(x)\right]dx=C_1C_3\gamma_nm_n^{1/2}\frac{2}{4-(2\alpha-a)\xi}+o\left(\gamma_nm_n^{1/2}\right).
\label{eq:24}
\end{eqnarray}
\item [$\bullet$] In the case  {$a\in \left(0,\frac{2}{9}\alpha\right)$}, we have $\gamma_nm_n^{1/2}=o\left(m_n^{-4}\right)$ and $\lim_{n\to+\infty}\left(n\gamma_n\right)>2a$.{Then the application of} lemma \ref{lem:1} gives
\begin{eqnarray}
\int^{1}_{0}Var\left[f_n(x)\right]dx=\Pi_n^2\displaystyle\sum_{k=1}^{n}\Pi_k^{-2}o\left(m^{-4}_{k}\right)=o\left(m^{-4}_{n}\right).
\label{eq:25}
\end{eqnarray}
\end{itemize}
Part 1 of Proposition \ref{prop:3} follows from the combination of (\ref{eq:18}) and (\ref{eq:25}), Part 2 from that of (\ref{eq:18}) and (\ref{eq:24}), Part 3 from (\ref{eq:19}) and (\ref{eq:24}).
\subsection{ Proof Corollary \ref{coro:1}}\label{sub:6.3}
Set
\begin{eqnarray*}
K_1(\xi)=\frac {C_1C_3}{2-(\alpha-\frac{a}{2})\xi},\quad K_2(\xi)=\frac{4C_2}{(1-2\alpha\xi)^2}.
\end{eqnarray*}
It follows from Proposition \ref{prop:3} that
\begin{eqnarray}
MISE\left[f_n\right]= 
\begin{cases}
K_2(\xi)m_n^{-4}\left[1+o(1)\right]&\quad\text{if}\,{a\in \left(0,\frac{2}{9}\alpha\right)}\\  m_n^{-4}\left[K_2(\xi)+K_1(\xi)\g_nm_n^{9/2}+o(1)\right]&\quad\text{if}\,a=\frac{2}{9}\alpha\\
K_1(\xi)\g_nm_n^{1/2}\left[1+o(1)\right]&\quad\text{if}\, {a\in \left(\frac{2}{9}\alpha,1\right)}.
\end{cases}
\label{minMISE}
\end{eqnarray}
If $a=\frac{2}{9}\alpha$, $\left(K_2(\xi)m_n^{-4}+K_1(\xi)\g_nm_n^{1/2}\right)\in\mathcal{GS}\left(-\frac{8}{9}\alpha\right)$ with $\frac{8}{9}\alpha=4a$. If ${a\in \left(0,\frac{2}{9}\alpha\right)}$, $\left(K_2(\xi)m_n^{-4}\right)\in\mathcal{GS}(-4a)$ with $-4a>-\frac{8}{9}\alpha$, and, if ${a\in \left(\frac{2}{9}\alpha,1\right)}$, $\left(K_1(\xi)\g_nm_n^{1/2}\right)\in\mathcal{GS}\left(-(2\alpha-a)/2\right)$ with $-(2\alpha-a)/2>-4a$.
 It follows that, for a given $\alpha$, to minimize the $MISE\left[f_n\right]$, the parameter $a$ must be chosen equal to $\frac{2}{9}\alpha$. Moreover, in view of (\ref{minMISE}) the parameter $\alpha$ must be equal to 1. We conclude that to minimize the $MISE\left[f_n\right]$, the stepsize $(\g_n)$ must be chosen in $\mathcal{GS}(-1)$ and the order $(m_n)$ in $\mathcal{GS}(2/9)$.\\
Since the function $
x\mapsto K_2(\xi)x^{-4}+K_1(\xi)\g_nx^{1/2}$ attains its minimum at the point $\left[\frac{8K_2(\xi)}{\g_nK_1(\xi)}\right]^{2/9}$, to minimize $MISE\left[f_n\right]$, the order $(m_n)$ must be equal to $\left(\frac{8K_2(\xi)}{K_1(\xi)}\right)^{2/9}\g_n^{-2/9}$. For such a choice, the $MISE$ of $f_n$ becomes 
\begin{eqnarray*}
MISE\left[f_n\right]=\frac{9}{8}8^{1/9}K_1(\xi)^{8/9}K_2(\xi)^{1/9}\g_n^{8/9}[1+o(1)].
\end{eqnarray*}
Now, we assume that $\left(\gamma_n\right)=\left(\gamma_0n^{-1}\right)$ for some $\gamma_0\in \left(0,\infty\right)$. In this case, note that $\xi=\lim_{n\to \infty}\left(n\gamma_n\right)^{-1}=\gamma_0^{-1}$. Then, (\ref{eq:30}) (hence (\ref{eq:31})) is a consequence of the standard trade-off argument. 

\subsection{Proof Theorem \ref{theo:1}}
To prove theorem \ref{theo:1}, we will use the fact that if $x\in(0,1)$ and {$a\in \left[\frac{2}{9}\alpha,1\right]$}, then
\begin{eqnarray}
\gamma^{-1/2}_{n}m_n^{-1/4}(f_{n}(x)-\mathbb{E}[f_{n}(x)])\stackrel{\mathcal{D}}{\rightarrow}\mathcal{N}\left(0,\frac{2}{4-(2\alpha-a)\xi}C_3f(x)\psi(x)\right),
\label{eq:29}
\end{eqnarray}
which will be proved later. In the case {$a\in \left(0,\frac{2}{9}\alpha\right)$}, we have
\begin{eqnarray}
\gamma_n^{-1/2}m_{n}^{-1/4}(f_{n}(x)-f(x))=\gamma^{-\frac{1}{2}}_{n}m_n^{-1/4}(f_{n}(x)-\mathbb{E}[f_{n}(x)])+\gamma^{-1/2}_{n}m_n^{-1/4}(\mathbb{E}[f_{n}(x)]-f(x))\nonumber\\
=\gamma^{-1/2}_{n}m_n^{-1/4}(f_{n}(x)-\mathbb{E}[f_{n}(x)])-\gamma^{-1/2}_{n}m^{-9/4}_{n}\frac{2}{1-2a\xi}[\Delta_2(x)+o(1)],\nonumber
\end{eqnarray}
{then} if $\gamma^{-1/2}_{n}m^{-9/4}_{n}\rightarrow c$, for some $c\geq0$, Part 1 of theorem \ref{theo:1} follows from (\ref{eq:29}). In the case when $a=\frac{2}{9}\alpha$, Parts 1 and 2 of Theorem \ref{theo:1} follow from the combination of (\ref{eq:2}) and (\ref{eq:29}). In the case {$a\in \left(\frac{2}{9}\alpha,1\right]$}, Part 1 of Theorem \ref{theo:1} follows from the combination of (\ref{eq:3}) and (\ref{eq:29}).\\
Now in the case $\gamma_n^{-1/2}m_{n}^{-1/4}\rightarrow\infty$ and {$a\in \left(0,\frac{2}{9}\alpha\right)$}, we have
\begin{eqnarray}
m_{n}^{2}(f_{n}(x)-f(x))&=&m_{n}^{2}(f_{n}(x)-\mathbb{E}[f_{n}(x)])+m_{n}^2(\mathbb{E}[f_{n}(x)]-f(x))\nonumber\\
&=&m_{n}^{2}(f_{n}(x)-\mathbb{E}[f_{n}(x)])-\frac{2}{1-a\xi}\Delta_2(x)[1+o(1)].\nonumber
\end{eqnarray}
Noting that the equation (\ref{eq:5}) implies 
\begin{eqnarray*}
m_n^{2}\left(f_n(x)-\mathbb{E}\left[f_n(x)\right]\right)\stackrel{\mathbb{P}}{\rightarrow}0,
\end{eqnarray*}
then, we obtain Part 2 of Theorem \ref{theo:1}.
We now prove (\ref{eq:29}). We have
\begin{eqnarray}
f_{n}(x)-\mathbb{E}[f_{n}(x)]&=&(1-\gamma_{n})(f_{n-1}(x)-\mathbb{E}[f_{n-1}(x)])+\gamma_{n}(Z_{n}(x)-\mathbb{E}[Z_{n}(x)])\nonumber\\
&=&\Pi_{n}\displaystyle\sum^{n}_{k=1}\Pi^{-1}_{k}\gamma_{k}(Z_{k}(x)-\mathbb{E}[Z_{k}(x)]).\nonumber
\end{eqnarray}
We set
\begin{eqnarray*}
 Y_{k}(x)=\Pi^{-1}_{k}\gamma_{k}(Z_{k}(x)-\mathbb{E}[Z_{k}(x)]).
\end{eqnarray*}
The application of Lemma \ref{lem:1} ensures that
\begin{eqnarray*}
v_{n}^{2}&=&\sum^{n}_{k=1}Var[Y_{k}(x)] \nonumber \\
&=&\sum^{n}_{k=1}\Pi^{-2}_{k}\gamma_{k}^{2}Var[Z_{k}(x)] \nonumber \\
&=&\sum^{n}_{k=1}\Pi^{-2}_{k}\gamma_{k}^{2}m_k^{1/2}\left[C_3f(x)\psi(x)+o(1)\right]\\
&=&\frac{\gamma_{n}}{\Pi_{n}^{2}}m_n^{1/2}\left[\frac{2}{4-(2\alpha-a)\xi}C_3f(x)\psi(x)+o(1)\right].
\end{eqnarray*}
On the other hand, for all $p>0$, we have
\begin{eqnarray*}
\mathbb{E}[|Z_{k}(x)|^{2+p}]=O(m_n^{2+p}),
\end{eqnarray*}
and, since $\lim_{n\rightarrow\infty}(n\gamma_{n})>(\alpha-\frac{a}{2})/2=(2\alpha-a)/4$, there existe $p>0$ such that $\lim_{n\rightarrow\infty}(n\gamma_{n})>\frac{1+p}{2+p}(\alpha-\frac{a}{2})=\frac{(1+p)\alpha-\frac{1+p}{2}a}{p+2}>\frac{(1+p)\alpha-(2+p)a}{p+2}$, {then the application of} lemma \ref{lem:1} gives
\begin{eqnarray}
\sum^{n}_{k=1}\mathbb{E}[|Y_{k}(x)|^{2+p}]&=& O\left(\sum^{n}_{k=1}\Pi^{-2-p}_{k}\gamma_{k}^{2+p}\mathbb{E}[|Z_{k}(x)|^{2+p}]\right)\nonumber\\
&=&O\left(\sum^{n}_{k=1}\frac{\gamma_{k}^{2+p}}{\Pi^{2+p}_{k}}m_k^{2+p}\right)\nonumber\\
&=&O\left(\frac{\gamma_{n}^{1+p}}{\Pi^{2+p}_{n}}m_n^{2+p}\right),
\end{eqnarray}
hence
\begin{eqnarray}
\frac{1}{v_{n}^{2+p}}\sum^{n}_{k=1}\mathbb{E}[|Y_{k}(x)|^{2+p}]=O\left(\gamma^{p/2}_{n}m_n^{\frac{3}{4}(2+p)}\right).\nonumber
\end{eqnarray}
Then the convergence in (\ref{eq:29}) follows from the application of Lyapounov's Theorem.

\section*{Funding}
This work benefited from the financial support of the GDR 3477 GeoSto.

\section*{Acknowledgments}

The authors would like to thank the Editor and the referees for their very helpful comments, which led to considerable
improvement of the original version of the paper and a more sharply focused presentation.


\begin{thebibliography}{}

\bibitem[{Babu et al.(2002)}]{Bab02}
{Babu, G. J.} {Canty, A. J.}, and {Chaubey, Y. P.} (2002),
\newblock{'Application of Bernstein polynomials for smooth estimation of a distribution and density function'},
\newblock \textit{Journal of Statistical Planning and Inference}, {\bf 105}, 377--392.

\bibitem[{Bojanic and Seneta(1973)}]{Boj73}
{Bojanic, R.}, and {Seneta, E.} (1973),
\newblock{'A unified theory of regularly varying sequences'},
\newblock \textit{Mathematische Zeitschrift}, {\bf 134}, 91--106.

\bibitem[{Chen(1996)}]{Che96}
{Chen, S, X.} (1996),
\newblock{'Empirical likelihood confidence intervals for nonparametric density estimation'},
\newblock \textit{Biometrika}, {\bf 83}, 329--341.

\bibitem[{Chen(1999)}]{Che99}
{Chen, S, X.} (1999),
\newblock{'Beta kernel estimators for density functions'},
\newblock \textit{Computational Statistics and Data Analysis}, {\bf 31}, 131--145.

\bibitem[{Chen(2000)}]{Che00}
{Chen, S, X.} (2000),
\newblock{'Probability density function estimation using gamma kernels'},
\newblock \textit{Annals of the Institute of Statistical Mathematics}, {\bf 52}, 471--480.

\bibitem[{Ghosal(2001)}]{Gho01}
{Ghosal, S.} (2000),
\newblock{'Convergence rates for density estimation with Bernstein polynomials'},
\newblock \textit{The Annals of Statistics}, {\bf 29}, 1264--1280.

\bibitem[{Duflo(1997)}]{Duf97}
{Duflo, M.} (1997),
\newblock{'Random iterative models'},
\newblock \textit{Collection Applications of Mathematics. Springer, Berlin}.

\bibitem[{Galambos and Seneta(1973)}]{Gal73}
{Galambos, J.}, and {Seneta, E.} (1973),
\newblock{'Regularly varying sequences'},
\newblock \textit{Proceedings of the American Mathematical Society}, {\bf 41}, 110--116.

\bibitem[{H\"ardle(1991)}]{Har91}
{H\"ardle, W.} (1991),
\newblock{'Smoothing techniques with implementation in S'},
\newblock \textit{Springer, New York}.

\bibitem[{Hirukawa(2010)}]{Hir10}
{Hirukawa, M.} (2010),
\newblock{'Nonparametric multiplicative bias correction for kernel-type density  estimation on the unit Interval'},
\newblock \textit{Computational Statistics and Data Analysis}, {\bf 54}, 473--495.

\bibitem[{Igarashi and Kakizawa(2014)}]{Kak14}
{Igarashi, G.} and {Kakizawa, Y.} (2014),
\newblock{'On improving convergence rate of Bernstein polynomial density estimator'},
\newblock \textit{Journal of Nonparametric Statistics}, {\bf 26 }, 61--84.

\bibitem[{Jones(1993)}]{Jon93}
{Jones, M. C.} (1993),
\newblock{'Simple boundary correction for density estimation kernel'},
\newblock \textit{Statistics and Computing}, {\bf 3}, 135--146.

\bibitem[{Kakizawa(2004)}]{Kak04}
{Kakizawa, Y.} (2004),
\newblock{'Bernstein polynomial probability density estimation'},
\newblock \textit{Journal of Nonparametric Statistics}, {\bf 16 }, 709--729.

\bibitem[{Leblanc(2010)}]{Leb10}
{Leblanc, A.} (2010),
\newblock{'A bias-reduced approach to density estimation using Bernstein polynomials'},
\newblock \textit{Journal of Nonparametric Statistics }, {\bf  22}, 459--475.

\bibitem[{Lejeune and Sarda(1992)}]{Lej92}
{Lejeune, M.} and {Sarda, P.} (1992),
\newblock{'Smooth estimators of distribution and density functions'},
\newblock \textit{Computational Statistics and Data Analysis}, {\bf 14}, 457--471.

\bibitem[{Mokkadem et al.(2009)Mokkadem, Pelletier and Slaoui(2009)}]{Mok09} 
{Mokkadem, A.}, {Pelletier, M.} and {Slaoui, Y.} (2009),
\newblock{'The stochastic approximation method for the estimation of a multivariate probability density'},
\newblock \textit{Journal of Statistical Planning and Inference}, {\bf 139}, 2459--2478.

\bibitem[{M\"uller(1991)}]{Mul91}
{M\"uller}, {H.-G.} (1991),
\newblock{'Smooth optimum kernel estimators near endpoints'},
\newblock \textit{Biometrika}, {\bf 78}, 521--530,

\bibitem[{M\"uller(1993)}]{Mul93}
{M\"uller}, {H.-G.} (1993),
\newblock{'On the boundary kernel method for non-parametric curve estimation near endpoints'},
\newblock \textit{Scandinavian Journal of Statistics}, {\bf 20}, 313--328,

\bibitem[{M\"uller and Wand(1994)}]{Mul94}
{M\"uller}, {H.-G.} and {Wang}, {J.-L.} (1994),
\newblock{'Hazard rate estimation under random censoring with varying kernels and bandwidths'},
\newblock \textit{Biometrics}, {\bf 50}, 61--76, 

\bibitem[Parzen(1962)]{Par62}
{Parzen, E.} (1962),
\newblock{'On estimation of probability density and mode'},
\newblock \textit {Annals of Mathematical Statistics}, {\bf 33}, 1065--1070.

\bibitem[Politis and Romano(1995)]{Pol95}
{Politis, D.N.} and {Romano, J.P.}  (1995),
\newblock{'Bias-corrected nonparametric spectral estimation'},
\newblock \textit {Journal of Times Series Analysis}, {\bf 16}, 67--103.

\bibitem[{Rao(2005)}]{Rao05}
{Rao, B. L. S. P.} (2005),
\newblock{'Estimation of distribution and density functions by generalized Bernstein polynomials'},
\newblock \textit{Indian Journal of Pure and Applied Mathematics}, {\bf 36 }, 63--88.

\bibitem[{Robbins and Monro(1951)}]{Rob51}
{Robbins, H.}, and {Monro, S.} (1951),
\newblock{'A stochastic approximation method'},
\newblock \textit{Annals of Mathematical Statistics}, {\bf 22}, 400--407.

\bibitem[Rosenblatt(1956)]{Ros56}
{Rosenblatt, M.} (1956),
\newblock{'Remarks on some nonparametric estimates of density functions'},
\newblock \textit {Annals of Mathematical Statistics}, {\bf 27}, 832--837.

\bibitem[{Schucany, Gray and Owen(1971)}]{Sch71} 
{Schucany, W.R.}, {Gray, H.L.} and {Owen, B.D.} (1971),
\newblock{'On bias reduction in estimation'},
\newblock \textit {Journal of the American Statistical Association}, {\bf 66}, 524--533.

\bibitem[{Schucany and Sommers(1977)}]{Sch77} 
{Schucany, W.R.} and {Sommers, J.~P.} (1977),
\newblock{'Improvement of kernel type density estimators'},
\newblock \textit {Journal of the American Statistical Association}, {\bf 72}, 420--423.

\bibitem[{Schuster(1985)}]{Sch85} 
{Schuster, E.~F.} (1985),
\newblock{'Incorporating support constraints into nonparametric estimators of densities'},
\newblock \textit {Communications in Statistics. Theory and Methods }, {\bf 14}, 1123--1136.

\bibitem[{Silverman(1986)}]{Sil86} 
{Silverman, B. W.} (1986),
\newblock{'Density estimation for statistics and data analysis'},
\newblock \textit{Chapman and Hall}, London.


\bibitem[{Slaoui(2013)}]{Sla13} 
{Slaoui, Y.} (2013),
\newblock{'Large and moderate deviation principles for recursive kernel density estimators defined by stochastic approximation method'},
\newblock \textit{Serdica Matheatical Journal}, {\bf 39},~53--82.

\bibitem[{Slaoui(2014a)}]{Sla14a} 
{Slaoui, Y.} (2014a),
\newblock{'Bandwidth selection for recursive kernel density estimators defined by stochastic approximation method'},
\newblock \textit{Journal of Probability and Statistics}, {\bf 2014}, ID 739640, doi:10.1155/2014/739640.

\bibitem[{Slaoui(2014b)}]{Sla14b} 
{Slaoui, Y.} (2014b),
\newblock{'The stochastic approximation method for the estimation of a distribution function'}.
\newblock \textit{Mathematical Methods of Statistics}, {\bf 23},~306--325.


\bibitem[{Terrell and Scott(1980)}]{Ter80} 
{Terrell, G.R.} and {Scott, D.W.} (1980),
\newblock{'On improving convergence rates for nonnegative kernel density estimators'}.
\newblock \textit{The Annals of Statistics}, {\bf 8},~1160--1163.

\bibitem[{Vitale(1975)}]{Vit75}
{Vitale, R. A.} (1975),
\newblock{'A bernstein polynomial approach to density function estimation},
\newblock \textit{Statistical Inference and Related Topics'}, {\bf 2}, 87--99.

\end{thebibliography}
\end{document}